\documentclass[a4paper,11pt]{article}
\usepackage[top=2.5cm,bottom=2.5cm,left=2.2cm,right=2.2cm]{geometry}
\usepackage{amsfonts}
\usepackage{mathrsfs,amscd,amssymb,amsthm,amsmath,bm,graphicx,psfrag,subfigure,url,mathtools}
\usepackage{pict2e}
\usepackage{stfloats}
\usepackage{psfrag,amsmath}
\usepackage{tikz}
\usepackage{indentfirst}
\usepackage{hyperref}
\usepackage{bookmark}
\usepackage{enumerate}
\usepackage{latexsym,euscript,epic,eepic,color}
\usepackage{multirow}
\usepackage{multicol}
\usepackage{longtable}
\usepackage{adjustbox}
\usepackage[all]{xy}
\usepackage{setspace}
\usepackage{epstopdf}
\allowdisplaybreaks
\usepackage{authblk}
\usepackage{pifont}
\usepackage{mdframed}
\usepackage{booktabs,array}
\usepackage{tabularx}
\usepackage{hyperref}
\hypersetup{colorlinks=true, linkcolor=blue, filecolor=magenta, urlcolor=cyan,}
\usepackage[most]{tcolorbox}
\definecolor{mygray}{RGB}{240,240,240}
\tcbset{
  colback=mygray,
  boxrule=0pt,
}
\graphicspath{ {./images/} }

\newcommand{\customfootnote}[1]{
  \let\thefootnote\relax\footnotetext{#1}
}

\makeatletter

\renewcommand{\@seccntformat}[1]{{\csname the#1\endcsname}{\normalsize .}\hspace{.5em}}
\makeatother

\renewcommand{\thefootnote}{\fnsymbol{footnote}}
\usepackage{ifpdf}
\usepackage{indentfirst}

\def \[{\begin{equation}}
\def \]{\end{equation}}

\newcommand{\ex}{{\rm ex}}
\newcommand{\Ex}{{\rm EX}}
\newcommand{\spex}{{\rm spex}}
\newcommand{\SPEX}{{\rm SPEX}}
\newtheorem{thm}{Theorem}[section]

\newtheorem{claim}{Claim}

\newtheorem{lem}[thm]{Lemma}

\newtheorem{pb}{Problem}

\newtheorem{conj}{Conjecture}

\begin{document}
\baselineskip=0.23in

\title{\bf Spectral Tur\'an-type problem in non-$r$-partite graphs: Forbidden generalized book graph $B_{r,k}$}

\author[,1]{Yuantian Yu\thanks{ytyumath@sina.com (Y. Yu)}}
\author[,2,3]{Shuchao Li\thanks{Corresponding author, li@ccnu.edu.cn (S. Li)}}
\affil[1]{School of Science, East China University of Technology, Nanchang 330013, China}
\affil[2]{School of Mathematics and Statistics, and Hubei Key Lab--Math. Sci.,\linebreak Central China Normal University, Wuhan 430079, China}
\affil[3]{Key Laboratory of Nonlinear Analysis \& Applications (Ministry of Education),\linebreak Central China Normal University, Wuhan 430079, China}
\date{\today}
\maketitle

\begin{abstract}
Given a graph $H$, a graph is said to be $H$-free if it does not contain $H$ as a subgraph. A graph is color-critical when it has an edge whose removal leads to a reduction in its chromatic number. Nikiforov (2009) put forward a spectral-based version of a result by Simonovits. It was shown that if $H$ is a color-critical graph with a chromatic number of \(r + 1\), then \(T_r(n)\), the $r$-partite Tur\'an graph of order $n$, is the unique $H$-free graph of order $n$ that attains the maximum spectral radius. When dealing with a color-critical graph $H$ having a chromatic number of \(r + 1\), it becomes an interesting task to identify non-$r$-partite $H$-free graphs with the maximum spectral radius. For a graph $H$ with a chromatic number of \(r + 1\), we use \(\text{spex}_{r + 1}(n, H)\) to represent the maximum spectral radius among non-$r$-partite $H$-free graphs of order $n$. The set of all non-$r$-partite $H$-free graphs of order $n$ that have a spectral radius of \(\text{spex}_{r + 1}(n, H)\) is denoted as \(\text{SPEX}_{r + 1}(n, H)\). For \(r\geq2\) and \(k\geq1\), we define \(B_{r,k}\) as the graph constructed by connecting each vertex of \(K_r\) to every vertex of an independent set with size $k$. We refer to \(B_{r,k}\) as a book graph (in the case of \(r = 2\)) or a generalized book graph (when \(r\geq3\)). It should be noted that \(B_{r,k}\) is a color-critical graph with a chromatic number of \(r + 1\). Lin, Ning, and Wu (2021) identified the unique extremal graph within \(\text{SPEX}_3(n, B_{2,1})\); Li and Peng (2023) determined the unique extremal graph in \(\text{SPEX}_{r + 1}(n, B_{r,1})\) for all \(r\geq2\). Quite recently, Liu and Miao (2025) specified the unique extremal graph in \(\text{SPEX}_3(n, B_{2,k})\) for all \(k\geq2\). 
Inspired by these remarkable results, this paper, relying on spectral stability theory, local structure characterization, along with the theory of characteristic equations and Rayleigh quotient equations, aims to determine the unique extremal graph in \(\text{SPEX}_{r + 1}(n, B_{r,k})\) for \(r\geq3\), \(k\geq1\), and sufficiently large $n$. This work partially addresses an open problem put forward in \cite{YL}.

\vskip 0.2cm
\noindent {\bf Keywords:}  Non-$r$-partite graph; Generalized book graph; Spectral radius; Extremal graph\vspace{2mm}

\noindent {\bf AMS Subject Classification:} 05C50; 05C75
\end{abstract}


\section{\normalsize Introduction}\label{s1}
In this paper, we focus on simple and finite graphs exclusively. Unless specified otherwise, we adhere to the conventional notation and terminology (as can be seen, for example, in Bollob\'as \cite{BB1998}, Godsil and Royle~\cite{GGR2001}).

For a graph \(G = (V(G), E(G))\), we denote the \textit{order} (number of vertices) and \textit{size} (number of edges) of $G$ by \(|V(G)|\) and \(e(G) := |E(G)|\), respectively. When there is no risk of confusion, we also use “size” to refer to the cardinality of a set. As is standard, let \(K_n\) and \(C_n\) represent the complete graph and cycle of order $n$, respectively. A simple complete $r$-partite graph with $n$ vertices, where the vertex-partition classes have sizes \(n_1, n_2, \ldots, n_r\) (satisfying \(n_1 + n_2 + \cdots + n_r = n\)), is denoted by \(K_{n_1, n_2, \ldots, n_r}\). The $r$-partite Tur\'an graph of order $n$, denoted \(T_r(n)\), is the complete $r$-partite graph \(K_{n_1, n_2, \ldots, n_r}\) such that \(\sum_{i = 1}^r n_i = n\) and \(|n_i - n_j| \leq 1\) for all \(i \neq j\).

Given a graph $H$, a graph is said to be $H$-free if it contains no subgraph isomorphic to $H$. The \textit{Tur\'{a}n-type problem} aims to: determine the maximum number of edges, denoted \(\text{ex}(n, H)\), in an $n$-vertex $H$-free graph; and characterize the set \(\text{EX}(n, H)\), which consists of all $n$-vertex $H$-free graphs with exactly \(\text{ex}(n, H)\) edges. Here, \(\text{ex}(n, H)\) is referred to as the \textit{Tur\'an number} of $H$, and each graph in \(\text{EX}(n, H)\) is called an \textit{extremal graph} for $H$. Research on Tur\'an-type problems has garnered significant attention and has emerged as one of the most compelling fundamental topics in extremal graph theory (for surveys, see \cite{FS2013,V2011})

In 1941, Tur\'an determined the Tur\'an number of $K_{r+1},$ extending the result of Mantel \cite{M1907}. 
\begin{thm}[\cite{T1941}]\label{thm1.1}
For positive integers $r,n,$ we have $\ex(n,K_{r+1})=e(T_r(n))$ and $\Ex(n,K_{r+1})=\{T_r(n)\}.$ 
\end{thm}

A graph is said to be \textit{properly coloured} if each vertex is assigned a color such that the two endpoints of every edge have distinct colors. The \textit{chromatic number} \(\chi(G)\) of a graph $G$ is defined as the minimum number $k$ for which $G$ can be properly colored using $k$ colors. A graph $G$ is called \textit{color-critical} if there exists an edge $e$ in $G$ such that the chromatic number of the graph obtained by removing $e$ (denoted \(G - e\)) satisfies \(\chi(G - e) < \chi(G)\). It is worth noting that the complete graph \(K_{r+1}\) is color-critical and has a chromatic number of \(r + 1\). Building on this, Simonovits extended Theorem \ref{thm1.1} to the context of color-critical graphs.
\begin{thm}[\cite{Simo1968,Simo1974}]\label{thm1.2}
For a positive integer $r$ and a graph $H,$ if $H$ is color-critical with $\chi(H)=r+1$, then for sufficiently large $n,$ we have $\ex(n,H)=e(T_r(n))$ and $\Ex(n,H)=\{T_r(n)\}.$ 
\end{thm}
Let \(H\) be a graph with chromatic number \( \chi(H) = r + 1 \). We define \( \text{ex}_{r+1}(n, H) \) as the maximum number of edges in an \( n \)-vertex non-\( r \)-partite \( H \)-free graph, and let \( \text{EX}_{r+1}(n, H) \) denote the set of all \( n \)-vertex non-\( r \)-partite \( H \)-free graphs that contain exactly \( \text{ex}_{r+1}(n, H) \) edges.

By Theorem \ref{thm1.2}, for a color-critical graph \( H \) with chromatic number \( \chi(H) = r + 1 \), when \( n \) is sufficiently large, the unique extremal graph for \( H \) is the \( r \)-partite Tur\'an graph \( T_r(n) \). Inspired by this phenomenon, the authors of the current paper  \cite{YL} proposed the following problem.
\begin{pb}[\cite{YL}]\label{pb1}
Given a color-critical graph $H$ with $\chi(H)=r+1,$ determine the value of $\ex_{r+1}(n,H)$ and characterize all the extremal graphs in $\Ex_{r+1}(n, H).$
\end{pb}

Problem \ref{pb1} has garnered considerable attention in the research community. Relevant studies include those on \( K_3 \) (see \cite{CJ2002,E1955,AFGS,KP2005}), on odd cycles \( C_{2k+1} \) with \( k \geq 2 \) (see \cite{RWWY}), and on complete graphs \( K_{r+1} \) with \( r \geq 3 \) (see \cite{AFGS,KP2005}).

Let \( G \) be a graph with vertex set \( \{v_1, v_2, \ldots, v_n\} \). The adjacency matrix of \( G \), denoted \( A(G) \), is defined as an \( n \times n \) (0, 1)-matrix \( (a_{ij}) \) where \( a_{ij} = 1 \) if and only if \( v_i v_j \in E(G) \). It is evident that \( A(G) \) is a real symmetric matrix; consequently, its eigenvalues are real numbers and can be ordered in non-increasing sequence as $\lambda_1(G)\geq \lambda_2(G)\geq\cdots\geq\lambda_n(G).$ The \textit{spectral radius} of $G$ is defined by
$$
\rho(G)=\max\{|\lambda_i(G)|: 1\leq i\leq n\}.
$$
Note that \(A(G)\) is a non-negative matrix. By the Perron-Frobenius theorem, the \textit{spectral radius} \(\rho(G)\) equals the largest eigenvalue \(\lambda_1(G)\). Furthermore, if $G$ is connected, then \(A(G)\) is irreducible. In this case, there exists a unique unit positive eigenvector \(\boldsymbol{x} = (x_{v_1}, x_{v_2}, \ldots, x_{v_n})^T\) \linebreak corresponding to \(\rho(G)\), which we refer to as the \textit{Perron vector} of \(A(G)\).

The \textit{spectral Tur\'{a}n type problem} (also referred to as the \textit{Brualdi-Solheid-Tur\'{a}n type extremal problem}) aims to: determine the maximum spectral radius, denoted \(\text{spex}(n, H)\), among all $n$-vertex $H$-free graphs; and characterize the set \(\text{SPEX}(n, H)\), which consists of all $n$-vertex $H$-free graphs whose spectral radius equals \(\text{spex}(n, H)\). Here, \(\text{spex}(n, H)\) is called the \textit{spectral Tur\'an} number of $H$, and each graph in \(\text{SPEX}(n, H)\) is termed a \textit{spectral extremal graph} for $H$. Research on the spectral Tur\'an-type problem has received significant attention over the past fifteen years (for a survey, see \cite{LLF2022}).

In 2007, Nikiforov showed a spectral analog of the Tur\'an's theorem (Theorem \ref{thm1.1}).
\begin{thm}[\cite{N2007}]\label{thm1.3}
For positive integers $r,n,$ we have $\spex(n,K_{r+1})=\rho(T_r(n))$ and $\SPEX(n,$ $K_{r+1})=\{T_r(n)\}.$ 
\end{thm}

In 2009, Nikiforov proposed a spectral version of Simonovits' theorem (Theorem \ref{thm1.2}), thereby extending \ref{thm1.3} to the context of color-critical graphs.
\begin{thm}[\cite{Nik2009}]\label{thm1.4}
For a positive integer $r$ and a graph $H,$ if $H$ is color-critical with $\chi(H)=r+1$, then for sufficiently large $n,$ we have $\spex(n,H)=\rho(T_r(n))$ and $\SPEX(n,H)=\{T_r(n)\}.$ 
\end{thm}

Let \( H \) be a graph with chromatic number \( \chi(H) = r + 1 \). We denote by \( \text{spex}_{r+1}(n, H) \) the maximum spectral radius among all non-\( r \)-partite \( H \)-free graphs of order \( n \), and by \( \text{SPEX}_{r+1}(n, H) \) the set of all such non-\( r \)-partite \( H \)-free graphs of order \( n \) that have a spectral radius equal to \( \text{spex}_{r+1}(n, H) \).

By Theorem \ref{thm1.4}, for a color-critical graph \( H \) with chromatic number \( \chi(H) = r + 1 \), when \( n \) is sufficiently large, the unique spectral extremal graph for \( H \) is the \( r \)-partite Tur\'an graph \( T_r(n) \). Inspired by this elegant result, it is worthwhile to investigate the spectral counterpart of Problem \ref{pb1} put forward by Yu and Li \cite{YL}.

\begin{pb}[\cite{YL}]\label{pb2}
Given a color-critical graph $H$ with $\chi(H)=r+1,$ determine the value of $\spex_{r+1}(n,H)$ and characterize all the extremal graphs in $\SPEX_{r+1}(n, H).$
\end{pb}
Zhang and Zhao \cite{ZZ2023} solved Problem \ref{pb2} for odd cycles \( C_{2k+1} \) with \( k \geq 2 \). Let \( H_1 \) and \( H_2 \) be two graphs. We define \( H_1 \cup H_2 \) as their disjoint union, and \( H_1 \vee H_2 \) as their \textit{join}—this is constructed from \( H_1 \cup H_2 \) by connecting every vertex of \( H_1 \) to every vertex of \( H_2 \) with an edge. For positive integers \( r \geq 2 \) and \( k \geq 1 \), we define \( B_{r,k} := K_r \vee kK_1 \). When \( r = 2 \) (resp. \( r \geq 3 \)), \( B_{r,k} \) is referred to as a \textit{book graph} (resp. \textit{generalized book graph}). It is noteworthy that \( B_{r,k} \) is a color-critical graph with chromatic number \( r+1 \). The authors of  \cite{YL} solved Problem \ref{pb1} for generalized book graphs \( B_{r,k} \) with \( r \geq 3 \) and \( k \geq 1 \); Lin, Ning, and Wu \cite{LNW} solved Problem \ref{pb2} for \( B_{2,1} \); Li and Peng \cite{LP2023} solved Problem \ref{pb2} for \( B_{r,1} \) with \( r \geq 3 \); and recently, Liu and Miao \cite{LM2025} solved Problem \ref{pb2} for book graphs \( B_{2,k} \) with \( k \geq 1 \). In this paper, we focus on the spectral version of the main result in \cite{YL} and solve Problem \ref{pb2} for generalized book graphs \( B_{r,k} \) where \( r \geq 3 \) and \( k \geq 1 \).

Choose two partite sets \( T_1 \) and \( T_2 \) of the Turán graph \( T_r(n) \), where \( |T_1| = \left\lfloor \frac{n}{r} \right\rfloor \) and \( |T_2| = \left\lceil \frac{n}{r} \right\rceil \). We add a new edge \( uw \) within partite set \( T_2 \), and then remove \( \left\lfloor \frac{n}{r} \right\rfloor - 1 \) edges between vertex \( u \) and \( T_1 \), along with one edge between vertex \( w \) and \( T_1 \), such that in the resulting graph (denoted \( Y_r(n) \)), vertices \( u \) and \( w \) have no common neighbors in \( T_1 \). The main result of this paper is presented in the following.
\begin{thm}\label{thm1.5}
Let \( r \geq 3 \), \( k \geq 1 \), and let \( n \) be sufficiently large. If \( G \) is a non-\( r \)-partite \( B_{r,k} \)-free graph of order \( n \), then \( \rho(G) \leq \rho(Y_r(n)) \), with equality if and only if \( G = Y_r(n) \).
\end{thm}

\noindent{\bf Outline of the paper.}\ \ In the remainder of this section, we introduce the necessary notations and terminologies. Section \ref{s2} presents essential preliminary results. In Section \ref{s3}, we refine the structure of extremal graphs step by step and complete the proof of Theorem \ref{thm1.5}, leveraging methods from spectral stability theory, local structure characterization, as well as techniques involving characteristic equations and Rayleigh quotient equations. The final section contains concluding remarks.  \\

\noindent{\bf Notations and terminologies.}\ \ Let \( G \) be a graph. Two vertices \( u \) and \( v \) in \( G \) are said to be \textit{adjacent} (or \textit{neighbours}) if they are joined by an edge, denoted by \( u \sim v \). If \( uv \in E(G) \), let \( G - uv \) (resp. \( G - u \)) denote the graph obtained from \( G \) by deleting the edge \( uv \) (resp. the vertex \( u \)); this notation extends naturally to the deletion of multiple edges or vertices. Similarly, if \( uv \notin E(G) \), let \( G + uv \) denote the graph obtained from \( G \) by adding an edge between \( u \) and \( v \).  For two disjoint vertex subsets \( V_1 \) and \( V_2 \) of \( V(G) \), let \( G[V_1] \) denote the subgraph of \( G \) induced by \( V_1 \), and \( G[V_1, V_2] \) denote a bipartite subgraph of \( G \) induced by the edges between \( V_1 \) and \( V_2 \). The number of edges in \( G[V_1] \) and \( G[V_1, V_2] \) is abbreviated to \( e(V_1) \) and \( e(V_1, V_2) \), respectively. The set of \textit{neighbors} of a vertex \( v \) in \( G \) is denoted by \( N_G(v) \); its size, called the \textit{degree} of \( v \) in \( G \), is denoted by \( d_G(v) \). Let \( N_G[v] = N_G(v) \cup \{v\} \) (the closed neighborhood of \( v \)). For a vertex \( v \in V(G) \) and a subset \( W \subseteq V(G) \), denote \( N_W(v) = N_G(v) \cap W \) and \( d_W(v) = |N_W(v)| \) (the degree of \( v \) within \( W \)).  For a positive integer \( t \), the set \( \{1, 2, \ldots, t\} \) is abbreviated to \( [t] \).  
\section{\normalsize Preliminaries}\label{s2}
In this section, we present several fundamental lemmas and properties that will be instrumental in the proof of our main result. These include classic results from spectral graph theory, structural properties of Tur\'an graphs, and techniques involving characteristic equations and Rayleigh quotient equations, which lay the groundwork for the technical analysis in subsequent sections.

Nikiforov \cite{N2009} established a spectral version of Simonovits' stability theorem. As a corollary, Wang, Kang and Xue \cite{WKX2023} put forward the following spectral stability result.
\begin{lem}[Spectral stability \cite{WKX2023}]\label{lem2.1} 
Let \( H \) be a graph with \( \chi(H) = r + 1 \geq 3 \). For every \( \varepsilon > 0 \), there exists a constant \( \delta > 0 \) and an integer \( n_0 \) such that if \( G \) is an \( H \)-free graph on \( n \geq n_0 \) vertices with \( \rho(G) \geq \left(1 - \frac{1}{r} - \delta\right)n \), then \( G \) can be obtained from \( T_r(n) \) by adding and deleting at most \( \varepsilon n^2 \) edges.
\end{lem}
\begin{lem}[Rayleigh quotient equation]\label{lem2.2}
Let $G$ be a graph of order $n$, the spectral radius of $G$ satisfies
$$
\rho(G)=\max_{\boldsymbol{x}\in\mathbb{R}^n}\frac{\boldsymbol{x}^TA(G)\boldsymbol{x}}{\boldsymbol{x}^T\boldsymbol{x}}
=\max_{\boldsymbol{x}\in\mathbb{R}^n}\frac{2\sum_{uv\in E(G)}x_ux_v}{\boldsymbol{x}^T\boldsymbol{x}}.
$$
\end{lem}

In 1986, Wilf provided the first result regarding the spectral version of Tur\'an theorem. 
\begin{lem}[\cite{Wilf1986}]\label{lem2.3}
For $r\geq 1,$ if $G$ is a $K_{r+1}$-free graph of order $n$, then $\rho(G)\leq \left(1-\frac{1}{r}\right)n.$ 
\end{lem}

\begin{lem}[\cite{SD2020}]\label{lem2.4}
Let $v\in V(G)$ with $d_G(v)\geq 1.$ Then $\rho(G)\leq\sqrt{\rho^2(G-v)+2d_G(v)-1}.$ Equality holds if and only if either $G= K_n$ or $G= K_{1,n-1}$ with $d_G(v)=1.$
\end{lem}
Let \( M \) be an \( n \times n \) real symmetric matrix, and let \( \pi: V = V_1 \cup \cdots \cup V_\ell \) be a partition of \( V = \{1, 2, \ldots, n\} \). Corresponding to the partition \( \pi \), \( M \) can be partitioned into the following block matrix:  
$$
M = \begin{pmatrix} M_{11} & \cdots & M_{1\ell} \\ \vdots & \ddots & \vdots \\ M_{\ell 1} & \cdots & M_{\ell\ell} \end{pmatrix}
$$
The \textit{quotient matrix} of \( M \) with respect to \( \pi \) is the \( \ell \times \ell \) matrix \( B = (b_{ij}) \), where \( b_{ij} \) denotes the average row sum of the block \( M_{ij} \). The partition \( \pi \) is called \textit{equitable} if, for all \( i, j \in [\ell] \), each block \( M_{ij} \) has constant row sums.
\begin{lem}[\cite{G1993}]\label{lem2.5}
Let $M$ be a real symmetric matrix and let $B$ be a quotient matrix of $M$ with respect to an equitable partition. Then the eigenvalues of $B$ are also the eigenvalues of $M.$ Furthermore, if $M$ is nonnegative and irreducible, then $\lambda_1(M)=\lambda_1(B).$ 
\end{lem}
\begin{lem}[\cite{WXH2005}]\label{lem2.6}
Let $G$ be a connected graph with $u,v\in V(G).$ Assume that $\boldsymbol{x}$ is the Perron vector of $A(G)$, and $S\subseteq N_G(v)\backslash N_G[u]$. Let $G'=G-\{vw:w\in S\}+\{uw:w\in S\}$. If $S\neq \emptyset$ and $\boldsymbol{x}_u\geq \boldsymbol{x}_v,$ then
$\rho(G')>\rho(G).$
\end{lem}
\begin{lem}[\cite{CFTZ}]\label{lem2.02}
Let $V_1,V_2,\ldots,V_t$ be $t$ finite sets. Then
$$
\left|\bigcap_{i=1}^tV_i\right|\geq\sum_{i=1}^t|V_i|-(t-1)\left|\bigcup_{i=1}^tV_i\right|.
$$
\end{lem}
\begin{lem}[\cite{LP2023}]\label{lem2.7}
Let $K_{n_1,n_2,\ldots,n_r}$ be the complete $r$-partite graph with parts $V_1,V_2,\ldots,V_r$ satisfying $\left|V_i\right|=n_i$ for every $i\in [r]$ and $\sum_{i=1}^rn_i=n-1.$ Let $G$ be an $n$-vertex graph obtained from $K_{n_1,n_2,\ldots,n_r}$ by adding a new vertex $u$, choosing $v\in V_1,\,w\in V_2,$ removing the edge $vw,$ and adding the edges $uv,uw,$ and $uz$ for every $z\in \bigcup_{i=3}^rV_i.$ Then $\rho(G)\leq\rho\left(Y_r(n)\right),$ equality holds if and only if $G=Y_r(n).$
\end{lem}
\begin{lem}\label{lem2.8}
The size of $Y_r(n)$ satisfies $e\left(Y_r(n)\right)\geq\left(1-\frac{1}{r}\right)\frac{n^2}{2}-\frac{n}{r}-\frac{r}{8}+1.$
\end{lem}
\begin{proof}
By a simple calculation (see \cite{WKX2023}), 
$$
e\left(T_r(n)\right)\geq\left(1-\frac{1}{r}\right)\frac{n^2}{2}-\frac{r}{8}.
$$ 
By the construction of $Y_r(n),$ 
$$
e\left(Y_r(n)\right)=e\left(T_r(n)\right)-\left\lfloor\frac{n}{r}\right\rfloor+1
\geq \left(1-\frac{1}{r}\right)\frac{n^2}{2}-\frac{n}{r}-\frac{r}{8}+1,
$$ 
as desired.
\end{proof}
\section{\normalsize Proof of Theorem~\ref{thm1.5}}\label{s3}
In this section, we present the proof of Theorem \ref{thm1.5} and determine the graph with the maximum spectral radius among all non-\( r \)-partite \( B_{r,k} \)-free graphs of order \( n \).

Fix integers \( r \), \( k \), \( n \) with \( r \geq 3 \), \( k \geq 1 \) and \( n \) being sufficiently large. Let \( G \) be a graph with the maximum spectral radius among all non-\( r \)-partite \( B_{r,k} \)-free graphs of order \( n \), and denote \( \rho := \rho(G) \). Our first lemma shows that \( G \) is connected.
\begin{lem}\label{lem3.1}
$G$ is connected.
\end{lem}
\begin{proof}
Suppose, to the contrary, that \( G \) is not connected. Let \( G_1 \) and \( G_2 \) be two components of \( G \) with \( \rho(G_1) = \rho \). Take a vertex \( u \in V(G_1) \) and a vertex \( v \in V(G_2) \), and construct \( G' = G + uv \). Then clearly, \( G' \) is a non-\( r \)-partite \( B_{r,k} \)-free graph. However, since \( G_1 \) is a proper subgraph of the connected graph \( G'[V(G_1) \cup V(G_2)] \), by the Perron-Frobenius theorem, we have \( \rho = \rho(G_1) < \rho(G'[V(G_1) \cup V(G_2)]) = \rho(G') \), which contradicts the choice of \( G \).
\end{proof}

Since \( G \) is connected, by the Perron-Frobenius theorem, we may let \( \boldsymbol{x} \) be the positive eigenvector of \( A(G) \) corresponding to \( \rho \) with \( \max \{ x_v : v \in V(G) \} = 1 \). Let \( z \) be a vertex in \( V(G) \) with \( x_z = 1 \). Our second lemma provides a bound for \( \rho \).
\begin{lem}\label{lem3.2}
It holds
$$
\rho >\left\{
        \begin{array}{ll}
          \frac{2}{3}n - \frac{7}{12}, & \text{if $r=3$}; \\
          \frac{r-1}{r}n - \frac{2}{r}-\frac{r}{4n}, & \text{if $r\geq4$}.
        \end{array}
      \right.
$$
\end{lem}
\begin{proof} 
We first consider $r=3.$ Since $Y_3(n)$ is non-3-bipartite $B_{3,k}$-free, by the choice of $G$ one gets $\rho\geq\rho\left(Y_3(n)\right).$ Let $V_1,\,V_2$ and $V_3$ be parts of $T_3(n),$ where $|V_1|=\left\lfloor\frac{n}{3}\right\rfloor, \,|V_2|=\left\lceil\frac{n}{3}\right\rceil.$ Take $v\in V_1$ and $u,w\in V_2,$ let $Y_3(n)=T_3(n)+uw-vw-\left\{w'u:w'\in V_1\backslash \{v\}\right\}$. Then the quotient matrix of $A\left(Y_3(n)\right)$ with respect to the equitable partition $$
\pi: V\left(Y_3(n)\right)=\{v\}\cup\{u\}\cup\{w\}\cup\left(V_1\backslash \{v\}\right)\cup\left(V_2\backslash \{u,w\}\right)\cup V_3
$$
is given by
$$
B=\left(
    \begin{array}{cccccc}
      0 & 1 & 0 & 0 & \left\lceil\frac{n}{3}\right\rceil-2 & n-\left\lfloor\frac{n}{3}\right\rfloor-\left\lceil\frac{n}{3}\right\rceil \\[3pt]
      1 & 0 & 1 & 0 & 0 & n-\left\lfloor\frac{n}{3}\right\rfloor-\left\lceil\frac{n}{3}\right\rceil \\[3pt]
      0 & 1 & 0 & \left\lfloor\frac{n}{3}\right\rfloor-1 & 0 & n-\left\lfloor\frac{n}{3}\right\rfloor-\left\lceil\frac{n}{3}\right\rceil \\[3pt]
      0 & 0 & 1 & 0 & \left\lceil\frac{n}{3}\right\rceil-2 & n-\left\lfloor\frac{n}{3}\right\rfloor-\left\lceil\frac{n}{3}\right\rceil \\[3pt]
      1 & 0 & 0 & \left\lfloor\frac{n}{3}\right\rfloor-1 & 0 & n-\left\lfloor\frac{n}{3}\right\rfloor-\left\lceil\frac{n}{3}\right\rceil \\[3pt]
      1 & 1 & 1 & \left\lfloor\frac{n}{3}\right\rfloor-1 & \left\lceil\frac{n}{3}\right\rceil-2 & 0 
    \end{array}
  \right).
$$

If $n\equiv 0 \pmod{3},$ then the characteristic polynomial of $B$ is
\begin{align*}
f_1(x)&=\frac{1}{729}\big[729 x^6+\left(-243 n^2+243 n-729\right) x^4+\left(-54 n^3+162 n^2-486 n\right) x^3\\
&\,\,\,\,\,\,\,+\left(27 n^3+324 n^2-2673 n+2187\right) x^2+\left(162 n^3-1296 n^2+2916 n-2916\right) x \\
&\,\,\,\,\,\,\,-135 n^3+1215 n^2-2430 n\big].
\end{align*}
By Lemma \ref{lem2.5}, $\rho\left(Y_3(n)\right)$ is the largest zero of $f_1(x).$ On the other hand, by a direct computation, 
$$
f_1\left(\frac{2}{3}n - \frac{7}{12}\right)=\frac{-73728 n^5+519168 n^4-8313728 n^3+26451216 n^2-34291872 n+9787393}{2985984}<0
$$
for sufficiently large $n$. Therefore, $\rho\geq\rho\left(Y_3(n)\right)>\frac{2}{3}n - \frac{7}{12}.$

If $n\equiv 1 \pmod{3},$ then the characteristic polynomial of $B$ is
\begin{align*}
f_2(x)&=\frac{1}{729}\big[729 x^6+\left(-243 n^2+243 n-729\right) x^4+\left(-54 n^3+162 n^2-648 n+540\right) x^3\\
&\,\,\,\,\,\,\,+\left(27 n^3+324 n^2-2430 n+2808\right) x^2+\left(162 n^3-1296 n^2+3564 n-3888\right) x \\
&\,\,\,\,\,\,\,-135 n^3+1215 n^2-3240 n+2160\big].
\end{align*}
By Lemma \ref{lem2.5}, $\rho\left(Y_3(n)\right)$ is the largest zero of $f_2(x).$ On the other hand, by a direct computation, 
$$
f_2\left(\frac{2}{3}n - \frac{7}{12}\right)=-\frac{\left(96 n^2-316 n+353\right) \left(768 n^3-832 n^2+64228 n-60577\right)}{2985984}<0
$$
for sufficiently large $n$. Therefore, $\rho\geq\rho\left(Y_3(n)\right)>\frac{2}{3}n - \frac{7}{12}.$

If $n\equiv 2 \pmod{3},$ then the characteristic polynomial of $B$ is
\begin{align*}
f_3(x)&=\frac{1}{729}\big[729 x^6+\left(-243 n^2+243 n-972\right) x^4+\left(-54 n^3+162 n^2-486 n-702\right) x^3\\
&\,\,\,\,\,\,\,+\left(27 n^3+324 n^2-2835 n+2700\right) x^2+\left(162 n^3-1296 n^2+2754 n-1620\right) x \\
&\,\,\,\,\,\,\,-135 n^3+1215 n^2-2025 n-3375\big].
\end{align*}
By Lemma \ref{lem2.5}, $\rho\left(Y_3(n)\right)$ is the largest zero of $f_3(x).$ On the other hand, by a direct computation, 
$$
f_3\left(\frac{2}{3}n - \frac{7}{12}\right)=-\frac{9 (32 n^2 - 116 n +143 ) (256 n^3- 192 n^2 + 28620 n  +4633 )}{2985984}<0
$$
for sufficiently large $n$. Therefore, $\rho\geq\rho\left(Y_3(n)\right)>\frac{2}{3}n - \frac{7}{12}.$

Now we consider $r\geq 4.$ Let \(H\) be a non-$r$-partite \({B}_{r,k}\)-free graph of order \(n\) with maximum number of edges. Then since $Y_r(n)$ is non-$r$-partite $B_{r,k}$-free. By the choice of $H$ and Lemma \ref{lem2.8}, one has \(e\left( H\right) \geq  e\left(Y_r(n)\right)\geq\left(1-\frac{1}{r}\right)\frac{n^2}{2}-\frac{n}{r}-\frac{r}{8}+1.\)

Together with the choice of $G$ and Lemma \ref{lem2.2}, we have

\begin{align*}
\rho  \geq  \rho \left( H\right)  \geq  \frac{{\bf 1}^{T}A\left( H\right) {\bf 1}}{{\bf 1}^{T}{\bf 1}} = \frac{{2e}\left( H\right) }{n} &\geq  \frac{2}{n}\left[\left(1-\frac{1}{r}\right)\frac{n^2}{2}-\frac{n}{r}-\frac{r}{8}+1\right] \\
&> \left( {1 - \frac{1}{r}}\right) n - \frac{2}{r} - \frac{r}{4n}.
\end{align*}
This completes the proof of this lemma.
\end{proof}

In the remainder of this section, let \(\varepsilon\) be a fixed constant with \(0 < \varepsilon  < \frac{1}{36r^8}\).

\begin{lem}\label{lem3.3} 
It holds that
$$
e\left( G\right)  \geq  e\left( {{T}_{r}\left( n\right) }\right)- \frac{\varepsilon}{2}{n}^{2},
$$
and that \(G\) admits a partition \(V\left( G\right) = {V}_{1} \cup  V_2\cup \cdots  \cup  {V}_{r}\) such that \(\sum_{1 \leq  i < j \leq  r}e\left( {{V}_{i},{V}_{j}}\right)\) attains the maximum, \(\sum_{i =  1}^{r}e\left( {V}_{i}\right)  \leq  \frac{\varepsilon }{2}{n}^{2}\) and for each \(i \in  \left\lbrack  r\right\rbrack\),
$$
\left( {\frac{1}{r} - 2\sqrt{\varepsilon }}\right) n < \left| {V}_{i}\right|  < \left( {\frac{1}{r} + 2\sqrt{\varepsilon }}\right) n.
$$
\end{lem}
\begin{proof}
Note that \( n \) is sufficiently large. By Lemmas \ref{lem2.1} and \ref{lem3.2}, \( G \) can be obtained from \( T_r(n) \) by adding and deleting at most \( \frac{\varepsilon}{2} n^2 \) edges. Hence, \( e(G) \geq e(T_r(n)) - \frac{\varepsilon}{2} n^2 \), and there exists a partition \( V(G) = U_1 \cup U_2 \cup \cdots \cup U_r \) satisfying:  \vspace{2mm}

$\bullet$\ \ \ \ \( \sum_{i=1}^r e(U_i) \leq \frac{\varepsilon}{2} n^2 \),  \vspace{2mm}

$\bullet$\ \ \ \  \( \sum_{1 \leq i < j \leq r} e(U_i, U_j) \geq e(T_r(n)) - \frac{\varepsilon}{2} n^2 \),\vspace{2mm}  

$\bullet$\ \ \ \  \( \left\lfloor \frac{n}{r} \right\rfloor \leq |U_i| \leq \left\lceil \frac{n}{r} \right\rceil \) for each \( i \in [r] \).  \vspace{2mm}

Let \( V(G) = V_1 \cup V_2 \cup \cdots \cup V_r \) be a partition maximizing \( \sum_{1 \leq i < j \leq r} e(V_i, V_j) \). Then  
$$
\sum_{1 \leq i < j \leq r} e(V_i, V_j) \geq \sum_{1 \leq i < j \leq r} e(U_i, U_j) \geq e(T_r(n)) - \frac{\varepsilon}{2} n^2
$$  
and  
$$
\sum_{i=1}^r e(V_i) \leq \sum_{i=1}^r e(U_i) \leq \frac{\varepsilon}{2} n^2.
$$  

Without loss of generality, we assume \( \left| |V_1| - \frac{n}{r} \right| = \max_{j \in [r]} \left| |V_j| - \frac{n}{r} \right| \). Suppose, to the contrary, that \( \left| |V_1| - \frac{n}{r} \right| \geq 2\sqrt{\varepsilon}\,n \). Then,  
\begin{align}
e(G) &\leq \sum_{1 \leq i < j \leq r} |V_i| \cdot |V_j| + \sum_{i=1}^r e(V_i)\notag \\
&\leq |V_1| \cdot (n - |V_1|) + \sum_{2 \leq i < j \leq r} |V_i| \cdot |V_j| + \frac{\varepsilon}{2}\,n^2\notag \\
&= |V_1| \cdot (n - |V_1|) + \frac{1}{2} \left[ \left( \sum_{i=2}^r |V_i| \right)^2 - \sum_{i=2}^r |V_i|^2 \right] + \frac{\varepsilon}{2}\,n^2 \notag\\
&\leq |V_1| \cdot (n - |V_1|) + \frac{1}{2}(n - |V_1|)^2 - \frac{1}{2(r-1)}(n - |V_1|)^2 + \frac{\varepsilon}{2}\,n^2 \notag\\
&= \frac{r-1}{2r}\,n^2 - \frac{1}{2r(r-1)}\,n^2 + \frac{1}{r-1} |V_1| \cdot n - \frac{r}{2(r-1)} |V_1|^2 + \frac{\varepsilon}{2}\,n^2\notag\\
&= \frac{r - 1}{2r}{n}^{2} - \frac{r}{2\left( {r - 1}\right) }{\left( \frac{n}{r} - \left| {V}_{1}\right| \right) }^{2} + \frac{\varepsilon }{2}{n}^{2}\notag\\
&\leq  \frac{r - 1}{2r}{n}^{2} - \frac{r}{2\left( {r - 1}\right) } \cdot  {4\varepsilon }{n}^{2} + \frac{\varepsilon }{2}{n}^{2}\notag\\
&< \frac{r - 1}{2r}{n}^{2} - \frac{3}{2}\varepsilon {n}^{2}. \label{q2}
\end{align}

On the other hand,
$$
e\left( G\right)  \geq  e\left( {{T}_{r}\left( n\right) }\right)  - \frac{\varepsilon }{2}{n}^{2} \geq  \frac{r - 1}{2r}{n}^{2} - \frac{r}{8} - \frac{\varepsilon }{2}{n}^{2} > \frac{r - 1}{2r}{n}^{2} - \varepsilon {n}^{2}
$$
for large enough \(n\), a contradiction to \eqref{q2}.
\end{proof}

For convenience, in the remainder of this section, for $v\in V(G),$ we denote by $d(v):=d_G(v)$; for $i\in[r],$ denote 
$$
W_i:=\left\{{v \in  {V}_{i} :  {d}_{{V}_{i}}\left( v\right)  \geq  3\sqrt{\varepsilon }n}\right\},\,\,\,\,\,\,\,\,\,\,
L_i:=\left\{ {v \in  V_i:  d\left( v\right)  \leq  \left( {1 - \frac{1}{r} - 5\sqrt{\varepsilon}}\right) n}\right\}.
$$
Let $W :=\bigcup_{i = 1}^{r}{W}_{i}$ and $L :=\bigcup_{i = 1}^{r}{L}_{i}.$
\begin{lem}\label{lem3.4} 
\(\left| L\right|  \leq  \sqrt{\varepsilon }n\) and \(W\subseteq L.\)
\end{lem}
\begin{proof}
We first show \(\left| L\right|  \leq  \sqrt{\varepsilon }n\). Suppose, to the contrary, that \(\left| L\right|  > \sqrt{\varepsilon }n\). Then there is a subset \(L' \subseteq  L\) with \(\left| L'\right|  = \left\lfloor  {\sqrt{\varepsilon }n}\right\rfloor\). Therefore,
\begin{align*}
e\left( {{V \backslash  {L}'}}\right)&\geq  e\left( G\right)- \sum_{{v \in  {L}'}}d\left( v\right)\\
&\geq  e\left( {T}_r(n)\right)  - \frac{\varepsilon }{2}{n}^{2} - \sqrt{\varepsilon }n \cdot  \left( {1 - \frac{1}{r} - 5\sqrt{\varepsilon }}\right) n\\
&\geq  \left( {1 - \frac{1}{r}}\right) \frac{{n}^{2}}{2} - \frac{r}{8} - \left( {1 - \frac{1}{r} - \frac{9}{2}\sqrt{\varepsilon }}\right) \sqrt{\varepsilon }{n}^{2}\\
&= \left( {1 - \frac{1}{r}}\right)  \cdot  \frac{{\left( n - \left\lfloor  \sqrt{\varepsilon }n\right\rfloor  \right) }^{2}}{2} + \left( {1 - \frac{1}{r}}\right)  \cdot  n \cdot  \left\lfloor  {\sqrt{\varepsilon }n}\right\rfloor   - \left( {1 - \frac{1}{r}}\right)  \cdot  \frac{{\left\lfloor  \sqrt{\varepsilon }n\right\rfloor  }^{2}}{2}\\
&\,\,\,\,\,\,\,
- \frac{r}{8} - \left( {1 - \frac{1}{r} - \frac{9}{2}\sqrt{\varepsilon }}\right) \sqrt{\varepsilon }{n}^{2}\\
&\geq  \left( {1 - \frac{1}{r}}\right)  \cdot  \frac{{\left(n-\left\lfloor  \sqrt{\varepsilon }n\right\rfloor\right) }^{2}}{2} + \left( {1 - \frac{1}{r}}\right)  \cdot  n \cdot  \left( {\sqrt{\varepsilon }n - 1}\right)  - \left( {1 - \frac{1}{r}}\right)  \cdot  \frac{\varepsilon {n}^{2}}{2}\\
&\,\,\,\,\,\,\,- \frac{r}{8} - \left( {1 - \frac{1}{r}}\right) \sqrt{\varepsilon }{n}^{2} + \frac{9}{2}\varepsilon {n}^{2}\\
&= \left( {1 - \frac{1}{r}}\right)  \cdot  \frac{{\left( n - \left\lfloor  \sqrt{\varepsilon }n\right\rfloor\right) }^{2}}{2} - \left( {1 - \frac{1}{r}}\right) n + \left( {8 + \frac{1}{r}}\right) \frac{\varepsilon {n}^{2}}{2} - \frac{r}{8}\\
&> e\left( {T}_r\left(n - \left\lfloor  \sqrt{\varepsilon }n\right\rfloor\right)\right).
\end{align*}
Note that \({B}_{r,k}\) is a color-critical graph with $\chi({B}_{r,k})=r+1$, and $n$ is sufficiently large. By Theorem \ref{thm1.2}, $e\left( {T}_r\left(n - \left\lfloor  \sqrt{\varepsilon }n\right\rfloor\right)\right)  = \ex\left({n - \left\lfloor  \sqrt{\varepsilon }n\right\rfloor,{B}_{r,k}}\right).$ Hence, \(e\left( { {V \backslash  {L}'} }\right)  > \ex\left({n - \left\lfloor  \sqrt{\varepsilon }n\right\rfloor,{B}_{r,k}}\right)\), which implies that \(G\left\lbrack  {V \backslash   L'}\right\rbrack\) and thus $G$ contain a copy of \({B}_{r,k}\), a contradiction. 

Next, we prove that \(W \subseteq  L\). To this end, we first establish the following claim, which provides an upper bound for $|W|$.
\begin{claim}\label{c1}
\(\left| W\right|  \leq  \frac{1}{3}\sqrt{\varepsilon }n\).
\end{claim}
\begin{proof}[\bf Proof of Claim \ref{c1}] 
It follows from Lemma \ref{lem3.3} that \(\sum_{i = 1}^{r}e\left( {V}_{i}\right)  \leq  \frac{\varepsilon }{2}{n}^{2}\). On the other hand, 
$$
2e\left(V_i\right) =\sum_{u \in  {V}_{i}}{d}_{{V}_{i}}\left( u\right)  \geq \sum_{u \in  {W}_{i}}{d}_{{V}_{i}}\left( u\right)  \geq  \left| {W}_{i}\right|  \cdot  3\sqrt{\varepsilon }n.
$$
Thus,
$$
\frac{\varepsilon}{2}{n}^{2} \geq \sum_{i = 1}^{r}e\left( {V}_{i}\right)  \geq  \left| W\right|  \cdot  \frac{3\sqrt{\varepsilon}}{2}n.
$$
So we obtain $\left| W\right|  \leq  \frac{1}{3}\sqrt{\varepsilon}n$.
\end{proof}

Now suppose, to the contrary, that there exists a vertex \( u_0 \in W \setminus L \). Without loss of generality, let \( u_0 \in V_1 \). Since the partition \( V(G) = V_1 \cup V_2 \cup \cdots \cup V_r \) maximizes \( \sum_{1 \leq i < j \leq r} e(V_i, V_j) \), we have \( d_{V_1}(u_0) \leq d_{V_i}(u_0) \) for each \( i \in [r] \setminus \{1\} \). Thus, \( d(u_0) \geq r \cdot d_{V_1}(u_0) \), i.e., \( d_{V_1}(u_0) \leq \frac{1}{r} d(u_0) \). On the other hand, since \( u_0 \notin L \), we get \( d(u_0) > \left(1 - \frac{1}{r} - 5\sqrt{\varepsilon}\right) n \). Hence,  
\begin{align}
{d}_{{V}_{2}}\left( {u}_{0}\right)  &= d\left( {u}_{0}\right)  - {d}_{{V}_{1}}\left( {u}_{0}\right)  - \sum_{{i = 3}}^{r}{d}_{{V}_{i}}\left( {u}_{0}\right)\notag\\
&\geq  \left( {1 - \frac{1}{r}}\right) d\left( {u}_{0}\right)  - \sum_{{i = 3}}^{r}\left| {V}_{i}\right|\notag\\
&>  \left( {1 - \frac{1}{r}}\right) \left( {1 - \frac{1}{r} - 5\sqrt{\varepsilon }}\right) n - \left( {r - 2}\right) \left( {\frac{1}{r} + 2\sqrt{\varepsilon }}\right)n\notag\\
&= \left( {\frac{1}{{r}^{2}} - \left( {{2r} + 1 - \frac{5}{r}}\right) \sqrt{\varepsilon }}\right) n.\label{eq:3.3}
\end{align}
Since \( |W| \leq \frac{1}{3} \sqrt{\varepsilon}\,n \) and \( |L| \leq \sqrt{\varepsilon}\,n \), it follows that  
$$
d_{V_2 \setminus (W \cup L)}(u_0) > \frac{n}{r^2} - \left(2r + \frac{7}{3} - \frac{5}{r}\right) \sqrt{\varepsilon}\,n.
$$

On the other hand, note that \(u_0 \in  W,\) and so \(d_{V_1}(u_0)  \geq  3\sqrt{\varepsilon}n>\left|L\cup W\right|\). Therefore, we can choose a vertex \(u_1 \in  N_{V_1}(u_0)  \backslash (W \cup  L)\) and obtain
\begin{align}
{d}_{{V}_{2}}\left( {u}_{1}\right)  &= d\left( {u}_{1}\right)  - {d}_{{V}_{1}}\left( {u}_{1}\right)  - \sum_{{i = 3}}^{r}{d}_{{V}_{i}}\left( {u}_{1}\right)\notag\\
&> \left( {1 - \frac{1}{r} - 5\sqrt{\varepsilon }}\right) n - 3\sqrt{\varepsilon }n - \sum_{{i = 3}}^{r}\left|{V}_{i}\right|\notag\\
&> \left( {1 - \frac{1}{r} - 8\sqrt{\varepsilon }}\right) n - \left( {r - 2}\right) \left( {\frac{1}{r} + 2\sqrt{\varepsilon }}\right) n\notag\\
&= \frac{n}{r} - \left( {{2r} + 4}\right) \sqrt{\varepsilon }n. \label{eq:3.4}
\end{align}
Therefore, by Lemmas \ref{lem2.02}, \ref{lem3.3} and \eqref{eq:3.3}-\eqref{eq:3.4},
\begin{align*}
\left| {{N}_{{V}_{2}}\left( {u}_{0}\right)  \cap  {N}_{{V}_{2}}\left( {u}_{1}\right) }\right|  &\geq  \left| {{N}_{{V}_{2}}\left( {u}_{0}\right) }\right|  + \left| {{N}_{{V}_{2}}\left( {u}_{1}\right) }\right|  - \left| {{N}_{{V}_{2}}\left( {u}_{0}\right)  \cup  {N}_{{V}_{2}}\left( {u}_{1}\right) }\right|\\
&> \frac{n}{{r}^{2}} - \left( {{2r} + 1 - \frac{5}{r}}\right) \sqrt{\varepsilon }n + \frac{n}{r} - \left( {{2r} + 4}\right) \sqrt{\varepsilon }n - \left( {\frac{1}{r} + 2\sqrt{\varepsilon }}\right) n\\
&> \frac{n}{{r}^{2}} - \left( {{4r} + 7}\right) \sqrt{\varepsilon }n.
\end{align*}
Since \(\left| W\right|  \leq  \frac{1}{3}\sqrt{\varepsilon }n,\,\left| L\right|  \leq  \sqrt{\varepsilon }n\), one has
\begin{align*}
\left| \left( {{N}_{{V}_{2}}\left( {u}_{0}\right)  \cap  {N}_{{V}_{2}}\left( {u}_{1}\right) }\right) \setminus \left({W \cup  L}\right) \right| &>\frac{n}{{r}^{2}} - \left( {{4r} + 7}\right) \sqrt{\varepsilon }n - \frac{4}{3}\sqrt{\varepsilon }n\\
&=\frac{n}{{r}^{2}} - \left( {{4r} + \frac{25}{3}}\right) \sqrt{\varepsilon }n > 0.
\end{align*}
Hence, there is a vertex \({u}_{2}\) in \({V}_{2} \setminus  \left( {W \cup  L}\right)\) adjacent to both \({u}_{0}\) and \({u}_{1}\). For an integer \(s\) with \(2 \leq  s \leq  r - 1\), assume that for any \(1 \leq  i \leq  s\), there is a vertex \(u_i\in V_i \setminus  \left( {W\cup L}\right)\) such that \(\left\{  {{u}_{0},{u}_{1},\ldots,{u}_{s}}\right\}\) is a clique. We next consider the number of common neighbors of these vertices in \({V}_{s + 1} \backslash  \left( W\cup L\right)\). Similarly, by \eqref{eq:3.3} and \eqref{eq:3.4}, we have 
$$
{d}_{{V}_{s + 1}}\left( {u}_{0}\right) > \frac{n}{{r}^{2}} - \left( {{2r} + 1 - \frac{5}{r}}\right) \sqrt{\varepsilon}n,
$$
and for each $i\in [s]$,
$$
{d}_{{V}_{s + 1}}\left( {u}_{i}\right)  > \frac{n}{r} - \left( {{2r} + 4}\right) \sqrt{\varepsilon }n.
$$
Together with Lemmas \ref{lem2.02} and \ref{lem3.3}, one gets
\begin{align*}
\left|N_{V_{s + 1}}(u_0)\bigcap \left(\bigcap_{i = 1}^sN_{V_{s + 1}}(u_i)\right)\right|  &\geq  {d}_{{V}_{s + 1}}\left( {u}_{0}\right)  + \mathop{\sum }\limits_{{i = 1}}^{s}{d}_{{V}_{s + 1}}\left( {u}_{i}\right)  - s \cdot  \left| {V}_{s + 1}\right|\\
&>\frac{n}{{r}^{2}} - \left( {{2r} + 1 - \frac{5}{r}}\right) \sqrt{\varepsilon}n + s \cdot  \left( {\frac{n}{r} - \left( {{2r} + 4}\right) \sqrt{\varepsilon }n}\right)  \\
&\,\,\,\,\,\,\,-s \cdot  \left( {\frac{1}{r} + 2\sqrt{\varepsilon }}\right) n\\
&>\frac{n}{{r}^{2}} - \left( {{2sr} + {2r+1} + {6s}}\right) \sqrt{\varepsilon }n. 
\end{align*}
Hence 
$$\left|N_{V_{s + 1}}(u_0)  \bigcap  \left(\bigcap_{i = 1}^{s}N_{V_{s + 1}}(u_i)\right)\bigg\backslash \left({W \cup  L}\right) \right| >  \frac{n}{{r}^{2}} - \left( {{2sr} + {2r} + {6s} + \frac{7}{3}}\right) \sqrt{\varepsilon }n > k.
$$ 
That is to say, \({u}_{0},{u}_{1},\ldots ,{u}_{s}\) have at least \(k\) common neighbors in \({V}_{s + 1} \backslash  \left( {W \cup  L}\right)\). Thus, for each \(2\leq i\leq r-1\), there exists a vertex \(u_i \in  V_i\setminus (W \cup  L)\) such that \(\left\{  {{u}_{0},{u}_{1},\ldots ,{u}_{r - 1}}\right\}\) is a clique, and they have \(k\) common neighbors \({u}_{r,1},{u}_{r,2},\ldots ,{u}_{r,k}\) in \({V}_{r}\). Now \(G\left\lbrack  \left\{  {{u}_{0},{u}_{1},\ldots ,{u}_{r - 1},{u}_{r,1},{u}_{r,2},\ldots ,{u}_{r,k}}\right\}  \right\rbrack\) is isomorphic to \({B}_{r,k}\), a contradiction. Hence, $W\subseteq L$.
\end{proof}

\begin{lem}\label{lem3.5}
\(\chi\left( G\right)  \geq  r + 1\).
\end{lem}
\begin{proof}
Denote by \(s := \chi\left( G\right)\). Since \(G\) is non-$r$-partite, \(s \neq  r\). If \(s \leq  r - 1\), then $G$ is $K_r$-free. By Lemma \ref{lem2.3}, 
$$
\rho \leq  \left( {1 - \frac{1}{2}}\right) n < \frac{2}{3}n - \frac{7}{12}
$$
for $r=3,$ and
$$
\rho \leq  \left( {1 - \frac{1}{r-1}}\right) n <\frac{r-1}{r}n - \frac{2}{r}-\frac{r}{4n}
$$
for $r\geq 4,$ a contradiction to Lemma \ref{lem3.2}.

Therefore, \(\chi\left( G\right)  \geq  r + 1\).
\end{proof}
\begin{lem}\label{lem3.6} 
For each \(i \in  \left\lbrack  r\right\rbrack,\,e\left( {{V}_{i} \backslash  L}\right)  = 0\).
\end{lem}
\begin{proof}
Suppose, for contradiction, that there is an \({i}_{0} \in  \left\lbrack  r\right\rbrack\) such that
\(e\left(V_{i_0}\setminus L\right)  \geq  1\). Without loss of generality, we may assume that \(e\left(V_1\setminus L\right)  \geq  1\). Let
\({u}_{0}{u}_{1}\) be an edge in \(G\left\lbrack V_1\setminus  L\right\rbrack\). Now \({u}_{0},{u}_{1} \notin  L\), and so by Lemma \ref{lem3.4}, \({u}_{0},{u}_{1} \notin  W\).
Hence, $d\left( {u}_{0}\right)  > \left( {1 - \frac{1}{r} - 5\sqrt{\varepsilon}}\right) n$ and ${d}_{{V}_{1}}\left({u}_{0}\right)  < 3\sqrt{\varepsilon }n.$ Together with Lemma \ref{lem3.3}, one has
\begin{align}
{d}_{{V}_{2}}\left( {u}_{0}\right)  &= d\left( {u}_{0}\right)  - {d}_{{V}_{1}}\left( {u}_{0}\right)  - \mathop{\sum }\limits_{{i = 3}}^{r}{d}_{{V}_{i}}\left( {u}_{0}\right)\notag\\
&> \left( {1 - \frac{1}{r} - 5\sqrt{\varepsilon }}\right) n - 3\sqrt{\varepsilon }n - \left( {r - 2}\right) \left( {\frac{1}{r} + 2\sqrt{\varepsilon }}\right) n\notag\\
&= \left( {\frac{1}{r} - \left( {{2r} + 4}\right) \sqrt{\varepsilon }}\right) n.\label{q2.7}
\end{align}
Similarly, \({d}_{{V}_{2}}\left( {u}_{1}\right)  >  \left( {\frac{1}{r} - \left( {{2r} + 4}\right) \sqrt{\varepsilon }}\right) n.\) Together with Lemmas \ref{lem2.02} and \ref{lem3.3}, one has
\begin{align*}
\left| {{N}_{{V}_{2}}\left( {u}_{0}\right)  \cap  {N}_{{V}_{2}}\left( {u}_{1}\right) }\right|  &\geq  {d}_{{V}_{2}}\left( {u}_{0}\right)  + {d}_{{V}_{2}}\left( {u}_{1}\right)  - \left| {{N}_{{V}_{2}}\left( {u}_{0}\right)  \cup  {N}_{{V}_{2}}\left( {u}_{1}\right) }\right|\\
&> \left( {\frac{2}{r} - 2\left( {{2r} + 4}\right) \sqrt{\varepsilon }}\right) n - \left( {\frac{1}{r} + 2\sqrt{\varepsilon }}\right) n\\
&= \left( {\frac{1}{r} - \left( {{4r} + {10}}\right) \sqrt{\varepsilon}}\right)n.
\end{align*}
Note that $|L|\leq\sqrt{\varepsilon }n$ (see Lemma \ref{lem3.4}). Hence, 
$$\left| {\left( {{N}_{{V}_2}\left( {u}_{0}\right)  \cap  {N}_{{V}_2}\left( {u}_{1}\right) }\right)  \setminus  L}\right|  >  \left( {\frac{1}{r} - \left( {{4r} + {11}}\right) \sqrt{\varepsilon }}\right) n > 0.
$$ 
That is to say, there is a vertex
\({u}_{2}\) in \(V_2\setminus L\) adjacent to both \({u}_{0}\) and \({u}_{1}\). For an integer \(s\) with \(2 \leq  s \leq  r-1\), 
assume that for any \(1 \leq  i \leq  s\), there is a vertex \({u}_{i} \in  {V}_{i} \setminus  L\) such that \(\{{{u}_{0},{u}_{1}},\ldots, u_s\}\) forms a clique. 

We next consider the number of common neighbors of these vertices in \({V}_{s + 1} \setminus  L\). Similarly, by \eqref{q2.7} we get 
$$
{d}_{{V}_{s + 1}}\left( {u}_{i}\right)  > \left( {\frac{1}{r} - \left( {{2r} + 4}\right) \sqrt{\varepsilon }}\right) n
$$
for each $i\in \{ 0,1,\ldots, s\}$. Hence, by Lemmas \ref{lem2.02} and \ref{lem3.3},
\begin{align*}
\left|\bigcap_{{i = 0}}^{s}{N}_{{V}_{s + 1}}\left( {u}_{i}\right)\right|  &\geq  \sum_{{i = 0}}^{s}d_{V_{s + 1}}\left( {u}_{i}\right)  - s \cdot  \left| {V}_{s + 1}\right|
\\
&> \left( {s + 1}\right) \left( {\frac{1}{r} - \left( {{2r} + 4}\right) \sqrt{\varepsilon }}\right) n - s\left( {\frac{1}{r} + 2\sqrt{\varepsilon }}\right) n\\
&= \frac{n}{r} - \left( {{2rs} + {6s} + {2r} + 4}\right) \sqrt{\varepsilon }n.
\end{align*}
Therefore, \({u}_{0},{u}_{1},\ldots ,{u}_{s}\) have at least \(\left|\bigcap_{{i = 0}}^{s}{N}_{{V}_{s + 1}}\left( {u}_{i}\right) \right|  - \left| L\right|  > \frac{n}{r} - ({2rs} + {6s} + {2r}\) +5) \(\sqrt{\varepsilon }n \geq  k\) common neighbors in \({V}_{s + 1}\setminus L\). Thus, for each \(2\leq i\leq r - 1\), there exists \({u}_{i} \in  {V}_{i} \setminus L\) such that \(\left\{  {{u}_{0},{u}_{1},\ldots ,{u}_{r - 1}}\right\}\) forms a clique, and they have \(k\) common neighbors \({u}_{r,1},{u}_{r,2},\ldots ,{u}_{r,k}\) in \({V}_{r}\). Now \(G\left\lbrack  \left\{{{u}_{0},{u}_{1},\ldots ,{u}_{r-1},{u}_{r,1},{u}_{r,2},\ldots ,{u}_{r,k}}\right\}  \right\rbrack\) is isomorphic to \({B}_{r,k}\), a contradiction. 
\end{proof}

Note that by Lemma \ref{lem3.5}, \(\chi \left( G\right) \geq r+ 1\). It follows from Lemma \ref{lem3.6} that \(L \neq  \emptyset\). 

\begin{lem}\label{lem3.7} 
For each \(u\in V\left( G\right)\), if \(G - u\) is non-$r$-partite, then \(x_{u} \geq  1 - \frac{\sqrt{\varepsilon}n}{\rho}\).
\end{lem}
\begin{proof}
Suppose, to the contrary, that there exists a vertex \( u \in V(G) \) such that \( G - u \) is non-\( r \)-partite, yet \( x_u < 1 - \frac{\sqrt{\varepsilon}\,n}{\rho} \). Recall that \( x_z = \max \{ x_y : y \in V(G) \} = 1 \). Then it follows that $d(z)\geq\sum_{w\sim z}x_{w}=\rho x_{z} =\rho.$ Combining this with Lemma \ref{lem3.2}, we find that $z\notin L.$ Let \( G' = G - \{ uv : v \in N_G(u) \} + \{ uw : w \in N_G(z) \setminus L \} \). Then \( G' \) is a non-\( r \)-partite, \( B_{r,k} \)-free graph. Note that by Lemma \ref{lem3.4}, \( d_L(z) \leq |L| \leq \sqrt{\varepsilon}\,n \). Then
\begin{align*}
\rho = \rho x_{z} = \sum_{{w \sim  z}}x_{w} &=\sum_{{w \sim  z,w \in  L}}x_{w} + \sum_{{w \sim  z,w \notin  L}}x_{w} \\
&\leq  \left| L\right|  + \sum_{{w \sim  z,w \notin  L}}x_{w} \\
&\leq  \sqrt{\varepsilon}n +\sum_{{w \sim  z,w \notin  L}}x_{w}.
\end{align*}
Hence, \(\sum_{{w \sim  z,w \notin  L}}x_{w} \geq  \rho- \sqrt{\varepsilon}n\). By Lemma \ref{lem2.2}, we have
\begin{align*}
\rho \left( {G}'\right)- \rho &\geq  \frac{\boldsymbol{x}^{T}\left( {A\left( {G}'\right)  - A\left( G\right) }\right) \boldsymbol{x}}{\boldsymbol{x}^T\boldsymbol{x}}\\
&= \frac{2x_{u}}{\boldsymbol{x}^{T}\boldsymbol{x}}\left( \sum_{{w \sim  z,w \notin  L}}x_{w} - \sum_{{v \sim  u}}x_{v}\right) \\
&\geq\frac{2x_{u}}{\boldsymbol{x}^{T}\boldsymbol{x}}\left( {\rho -\sqrt{\varepsilon}n - \rho x_{u}}\right) \\
&>0,
\end{align*}
a contradiction to the choice of \(G\). 
\end{proof}
\begin{lem}\label{lem3.8}
For each \(v \in  L\), it holds \(x_{v} < 1 - \frac{4\sqrt{\varepsilon}n}{\rho}\).
\end{lem}
\begin{proof}
Note that \(L \neq  \emptyset\). Let \(v\) be a vertex in \(L\). Then by the definition of \(L\), one has \(d\left( v\right) \leq  \left( {1 - \frac{1}{r} - 5\sqrt{\varepsilon}}\right) n\), and so
$$
\rho x_{v} =\sum_{{w \sim  v}}x_{w} \leq d\left( v\right) \leq  \left( {1 - \frac{1}{r} - {5\sqrt{\varepsilon}}}\right) n.
$$
Together with Lemma \ref{lem3.2}, one gets
$$
x_{v} \leq  \frac{\left( {1 - \frac{1}{r} - 5\sqrt{\varepsilon}}\right) n}{\rho} = \frac{\left( {1 - \frac{1}{r}}\right) n - \frac{2}{r} - \frac{r}{4n}}{\rho} - \frac{5\sqrt{\varepsilon}n - \frac{2}{r} - \frac{r}{4n}}{\rho }
< 1 - \frac{4\sqrt{\varepsilon}n}{\rho},
$$
as desired.
\end{proof}
\begin{lem}\label{lem3.9}
For each \(u \in  V\left( G\right)\), if \(G - u\) is \(r\)-{partite}, then
$$
{d}\left( u\right)  > \left\{  \begin{matrix} \frac{n}{18}, & \text{ for }r = 3; \\  
\left( {1 - \frac{4}{r} + \frac{3}{{r}^{2}}}\right) n - \frac{r}{4}, & \text{ for }r \geq  4. \end{matrix}\right.
$$
\end{lem}
\begin{proof}
Since \(G - u\) is $r$-partite, it is $K_{r+1}$-free. By Lemma \ref{lem2.3}, \(\rho\left( {G - u}\right) \leq  \left( {1 - \frac{1}{r}}\right) \left( {n - 1}\right)\). Combine with Lemmas \ref{lem2.4} and \ref{lem3.2}, for $r=3$, one gets
\begin{align*}
2{d}\left( u\right)  &\geq  {\rho }^{2} - {\rho }^{2}\left( {G - u}\right)  + 1 \\
&> {\left( \frac{2}{3}n - \frac{7}{12}\right) }^{2} - \frac{4}{9}{\left( n - 1\right) }^{2} + 1 \\
&= \frac{n}{9} + \frac{129}{144}\\
&> \frac{n}{9},
\end{align*}
and so \({d}\left( u\right)  > \frac{n}{18}\).

Similarly, for $r\geq 4,$ one has 
\begin{align*}
2{d}\left( u\right) &\geq  {\rho }^{2} - {\rho }^{2}\left( {G - u}\right) + 1 \\
&>  {\left[ \left( 1 - \frac{1}{r}\right) n - \frac{2}{r} - \frac{r}{4n}\right] }^{2} - {\left( 1 - \frac{1}{r}\right) }^{2}{\left( n - 1\right) }^{2} + 1 \\
&= {\left( 1 - \frac{1}{r}\right) }^{2}{n}^{2} - 2\left( {1 - \frac{1}{r}}\right) \left( {\frac{2}{r} + \frac{r}{4n}}\right) n + {\left( \frac{2}{r} + \frac{r}{4n}\right) }^{2} - {\left( 1 - \frac{1}{r}\right) }^{2}\left( {{n}^{2} - {2n} + 1}\right)  + 1\\
&> \left( {2 - \frac{8}{r} + \frac{6}{{r}^{2}}}\right) n -\frac{r}{2},
\end{align*}
and so \({d}\left( u\right)  > \left( {1 - \frac{4}{r} + \frac{3}{{r}^{2}}}\right) n -\frac{r}{4}\).
\end{proof}
\begin{lem}\label{lem3.10}
Either \(\left| L\right| = 1,\) or \(\sum_{{i = 1}}^{r}e\left( V_{i}\right) = 1\).
\end{lem}
\begin{proof}
Note that $|L|\geq 1$ and \(\sum_{{i = 1}}^{r}e\left( {V}_{i}\right)  \geq  1\). It suffices to show \(\sum_{{i = 1}}^{r}e\left( {V}_{i}\right) = 1\) if \(\left| L\right| \geq 2\). In fact, if \(\left| L\right|\geq2\) but \(\sum_{{i = 1}}^{r}e\left( {V}_{i}\right)  \geq  2\), then there exists a vertex \(u \in  L\) such that \(G - u\) is non-$r$-partite. By Lemma \ref{lem3.7}, \(x_{u} \geq  1 - \frac{\sqrt{\varepsilon}n}{\rho}\). However, by Lemma \ref{lem3.8}, \(x_{u} < 1 - \frac{4\sqrt{\varepsilon}n}{\rho}\), a contradiction.
\end{proof}
\begin{lem}\label{lem3.11}
\(\sum_{{i = 1}}^{r}e\left( {V}_{i}\right)=1\).
\end{lem}
\begin{proof}
Since \(\chi\left( G\right)  \geq  r + 1\), one has \(\sum_{{i = 1}}^{r}e\left( {V}_{i}\right)  \geq  1\). Suppose to the contrary that \(\sum_{{i = 1}}^{r}e\left( {V}_{i}\right)  \geq  2\). Then by Lemma \ref{lem3.10}, \(\left| L\right|  = 1\). Let \(L = \{ u\}\). Without loss of generality, we may assume \(u \in  {V}_{1}\). Then by Lemma \ref{lem3.6}, 
$$
{d}_{{V}_{1}}\left( u\right) = e\left( {V}_{1}\right) =\sum_{{i = 1}}^{r}e\left( {V}_{i}\right) \geq  2.
$$ 
To complete the proof of this Lemma, we need the following claims. 
\begin{claim}\label{c3}
For each $s\in[ r]$ and each vertex \(v \in  {N}_{{V}_{s}}\left( u\right)\), It holds that \(N\left( v\right)=\bigcup_{\substack{i = 1\\ i \neq  s}}^{r}{V}_{i}\). 
\end{claim}\begin{proof}[\bf Proof of Claim \ref{c3}] Assume, for contradiction, that there exist \( s \in [r] \) and a vertex \( u_s \in N_{V_s}(u) \) such that at least one vertex in \( \bigcup_{\substack{i = 1 \\ i \neq s}}^r V_i \) is not adjacent to \( u_s \). Construct  
$$
G' = G - uu_s \;+\; \left\{ u_s v \,:\, v \in \left( \bigcup_{\substack{i = 1 \\ i \neq s}}^r V_i \right) \setminus N_G(u_s) \right\}.
$$ 
Clearly, \( G' \) is non-\( r \)-partite and \( B_{r,k} \)-free.  

On the other hand, take a vertex \( w \in \left( \bigcup_{\substack{i = 1 \\ i \neq s}}^r V_i \right) \setminus N_G(u_s) \). The graph \( G - w \) is non-\( r \)-partite. By Lemma 3.7, \( x_w \geq 1 - \frac{\sqrt{\varepsilon}\,n}{\rho} \). Moreover, since \( u \in L \), Lemma 3.8 implies \( x_u \leq 1 - \frac{4\sqrt{\varepsilon}\,n}{\rho} \). By Lemma 2.2, 
\begin{align*}
\rho \left( {G}'\right)  - \rho \geq  \frac{\boldsymbol{x}^T\left( {A\left( {G}'\right)  - A\left( G\right) }\right) \boldsymbol{x}}{\boldsymbol{x}^{T}\boldsymbol{x}} &\geq  \frac{2x_{{u}_{s}}}{\boldsymbol{x}^{T}\boldsymbol{x}}\left( {x_{w} - x_{u}}\right)\\
&\geq  \frac{2 x_{{u}_{s}}}{\boldsymbol{x}^{T}\boldsymbol{x}}\left( {1 - \frac{\sqrt{\varepsilon}n}{\rho}-1 + \frac{4\sqrt{\varepsilon}n}{\rho}}\right)\\
&> 0,
\end{align*}
a contradiction to the choice of \(G\).
\end{proof}
\begin{claim}\label{c4}
For each \(i \in  \left[r\right]\), \({d}_{{V}_{i}}\left( u\right)  \leq  k - 1\).
\end{claim}
\begin{proof}[\bf Proof of Claim \ref{c4}] Assume, for contradiction, that there exists \( s \in [r] \) such that \( d_{V_s}(u) \geq k \). Since the partition \( V(G) = V_1 \cup V_2 \cup \cdots \cup V_r \) maximizes \( \sum_{1 \leq i < j \leq r} e(V_i, V_j) \), it follows that \( d_{V_j}(u) \geq d_{V_1}(u) \geq 2 \) for each \( j \in \{2, 3, \ldots, r\} \).  

Now, take \( u_j \in N_{V_j}(u) \) for all \( j \in [r] \setminus \{s\} \), and let \( u_{s,1}, u_{s,2}, \ldots, u_{s,k} \in N_{V_s}(u) \). By Claim \ref{c3}, the subgraph  
$$
G\left[ \left\{ u, u_1, \ldots, u_{s-1}, u_{s+1}, \ldots, u_r, u_{s,1}, \ldots, u_{s,k} \right\} \right]
$$  
is isomorphic to \( B_{r,k} \), a contradiction.  
\end{proof}

From Claim \ref{c4}, we have \( d(u) = \sum_{i=1}^r d_{V_i}(u) \leq (k - 1)r \). On the other hand, since \( L = \{u\} \), Lemma \ref{lem3.6}  implies \( G - u \) is \( r \)-partite. Then Lemma \ref{lem3.9} yields a contradiction to \( d(u) \leq (k - 1)r \). Thus, \( \sum_{i=1}^r e(V_i) = 1 \).
\end{proof}

By Lemma \ref{lem3.11}, we may assume \(e\left( {V}_{1}\right)= 1\) and \(e\left( {V}_{i}\right) = 0\) for each \(i \in  \{ 2,\ldots ,r\}\). Let $uv$ be the unique edge in \(G\left[ {V}_{1}\right]\). 
\begin{lem}\label{lem3.12}
There is an \(i \in  \{2,3,\ldots ,r\}\) such that \({N}_{{V}_{i}}\left( u\right) \cap {N}_{{V}_{i}}\left( v\right)  = \emptyset\). 
\end{lem}
\begin{proof}
Suppose to the contrary that for each \(i \in  \{ 2,3,\ldots ,r\}\), it holds \({N}_{{V}_{i}}\left( u\right)\cap {N}_{{V}_{i}}\left( v\right) \neq  \emptyset\). In order to induce a contradiction, we first consider the following two claims.

\begin{claim}\label{c5}
For each \(i \in \{2,3,\ldots,r\}\), \(\left| {{N}_{{V}_{i}}\left( u\right) \cap{N}_{{V}_{i}}\left( v\right) }\right| \leq  k - 1\). Consequently, \(\left| {{N}_{G}\left( u\right)  \cap  {N}_{G}\left( v\right) }\right|  \leq  \left( {r - 1}\right) \left( {k - 1}\right)\). 
\end{claim}
\begin{proof}[\bf Proof of Claim \ref{c5}]
Assume, for contradiction, that there exists \( s \in \{2, 3, \ldots, r\} \) such that \linebreak\( \left| N_{V_s}(u) \cap N_{V_s}(v) \right| \geq k \). Without loss of generality, let \( s = r \).  For each \( i \in \{2, 3, \ldots, r - 1\} \), fix \( u_i \in N_{V_i}(u) \cap N_{V_i}(v) \). Let \( u_{r,1}, u_{r,2}, \ldots, u_{r,k} \in N_{V_r}(u) \cap N_{V_r}(v) \).  

By an argument analogous to Claim 2, for each \( s \in \{2, 3, \ldots, r - 1\} \), \( N_G(u_s) = \bigcup_{\substack{i = 1 \\ i \neq s}}^r V_i \). Thus, the induced subgraph  
$
G\left[ \{u, v, u_2, \ldots, u_{r-1}, u_{r,1}, \ldots, u_{r,k}\} \right]
$
is isomorphic to \( B_{r,k} \), contradicting the \( B_{r,k} \)-freeness of \( G \).

It follows that for each \( i \in \{2, 3, \ldots, r\} \), \( \left| N_{V_i}(u) \cap N_{V_i}(v) \right| \leq k - 1 \). Summing over \( i = 2 \) to \( r \), we get:  
$$
\left| N_G(u) \cap N_G(v) \right| = \sum_{i=2}^r \left| N_{V_i}(u) \cap N_{V_i}(v) \right| \leq (r - 1)(k - 1),
$$ 
as desired.
\end{proof}
\begin{claim}\label{c6}
For each \(i \in  \{2,3,\ldots ,r\}\), it holds \({N}_{{V}_{i}}\left( u\right)  \cup  {N}_{{V}_{i}}\left( v\right)  = {V}_{i}\).
\end{claim}
\begin{proof}[\bf Proof of Claim \ref{c6}]
Suppose, for the sake of contradiction, that there exists some \( i \in \{2, 3, \ldots, r\} \) such that \( N_{V_i}(u) \cup N_{V_i}(v) \neq V_i \). Construct a new graph \( G' \) as follows:  
$$
G' = G \;+\; \left\{ vw \,:\, w \in V_i \setminus \bigl(N_{V_i}(u) \cup N_{V_i}(v)\bigr) \right\}.
$$
Evidently, \( G' \) remains non-\( r \)-partite and \( B_{r,k} \)-free.  

On the other hand, by the Perron-Frobenius theorem, the spectral radius satisfies \( \rho(G') > \rho(G) \). This contradicts the initial choice of \( G \) (as \( G \) was assumed to maximize the spectral radius among such graphs). 
\end{proof}

Now we come back to show Lemma \ref{lem3.12}.  By Lemma \ref{lem3.3} and Claim \ref{c5}, 
\begin{align}\notag
d\left( u\right) + {d}\left( v\right) &= \left| {{N}_{G}\left( u\right) \cap {N}_{G}\left( v\right) }\right| + \left| {{N}_{G}\left( u\right)}\cup {N}_{G}\left( v\right)\right| \\ \notag
&\leq  \left( {r - 1}\right) \left( {k - 1}\right) + 2 + \sum_{{i = 2}}^{r}\left| {V}_{i}\right| \\ \notag
&=\left( {r - 1}\right) \left( {k - 1}\right) + 2 + \left( {n - \left| {V}_{1}\right| }\right) \\ \label{eq:3.5}
&\leq \left( {1 - \frac{1}{r} + 2\sqrt{\varepsilon}}\right) n + \left( {r - 1}\right) \left( {k - 1}\right)  + 2.
\end{align}
By Lemma \ref{lem3.3} and Claim \ref{c6},
\begin{align}\label{eq:3.6}
{d}\left( u\right)  + {d}\left( v\right)  \geq  2 + \sum_{{i = 2}}^{r}\left| {V}_{i}\right| = 2 + n - \left| {V}_{1}\right|  \geq  \left( {1 - \frac{1}{r} - 2\sqrt{\varepsilon }}\right)n+2.
\end{align}

Without loss of generality, we may assume \({d}\left( u\right)  \leq  {d}\left( v\right)\). Note that \(G - u\) is \(r\)-partite, combine with Lemma \ref{lem3.9} and \eqref{eq:3.5}, for \(r = 3\), one has
\begin{align}\label{eq:3.7}
\frac{n}{18}<{d}\left( u\right)  \leq  \frac{1}{2}\left( {{d}\left( u\right)  + {d}\left( v\right) }\right)  \leq  \left( {\frac{1}{3} + \sqrt{\varepsilon }}\right) n + k,
\end{align}
for \(r \geq 4\), one has
\begin{align}\label{eq:3.8}
\left( {1 - \frac{4}{r} + \frac{3}{{r}^{2}}}\right) n - \frac{r}{4}<{d}\left( u\right)  \leq  \frac{1}{2}\left( {{d}\left( u\right)  + {d}\left( v\right) }\right)  \leq  \left(\frac{1}{2}-\frac{1}{2r}+\sqrt{\varepsilon }\right) n+\frac{1}{2}(r-1)(k-1)+1.
\end{align}
Combine with \eqref{eq:3.5}-\eqref{eq:3.8}, for \(r = 3\), one has
\begin{align}\label{eq:3.9}
\left( {\frac{1}{3} - \sqrt{\varepsilon }}\right) n + 1 \leq  {d}\left( v\right)  <  \left( {\frac{11}{18} + 2\sqrt{\varepsilon }}\right) n + {2k},
\end{align}
for $r\geq4,$ one has
\begin{align}\label{eq:3.10}
\left( {\frac{1}{2} - \frac{1}{2r} - \sqrt{\varepsilon }}\right) n + 1 \leq  {d}\left( v\right)  <  \left( {\frac{3}{r} - \frac{3}{{r}^{2}} + 2\sqrt{\varepsilon }}\right) n + {kr} - k + 3.
\end{align}

We complete the proof of this lemma by analyzing two cases.  

{\bf Case 1.} \(r = 3\). Suppose there exists \( i \in \{2, 3\} \) (say \( i = 2 \)) such that \( d_{V_2}(u) = 1 \). By \( N_{V_2}(u) \cap N_{V_2}(v) \neq \emptyset \) and Claim \ref{c6}, \( N_{V_2}(v) = V_2 \). Combining with Lemma \ref{lem3.3} and \eqref{eq:3.9}, we derive:  
$$
d_{V_3}(v) = d(v) - 1 - d_{V_2}(v) < \left( \frac{11}{18} + 2\sqrt{\varepsilon} \right)n + 2k - 1 - \left( \frac{1}{3} - 2\sqrt{\varepsilon} \right)n = \left( \frac{5}{18} + 4\sqrt{\varepsilon} \right)n + 2k - 1.
$$
Thus, the number of vertices in \( V_3 \) non-adjacent to \( v \) satisfies:  
$$
\left| V_3 \setminus N_{V_3}(v) \right| \geq \left( \frac{1}{3} - 2\sqrt{\varepsilon} \right)n - \left( \frac{5}{18} + 4\sqrt{\varepsilon} \right)n - 2k + 1 = \left( \frac{1}{18} - 6\sqrt{\varepsilon} \right)n - 2k + 1.
$$ 
Take \( w_1, w_2 \in V_3 \setminus N_{V_3}(v) \). By Claim \ref{c6}, \( w_1, w_2 \in N_{V_3}(u) \), so \( G - w_1 \) and \( G - w_2 \) are non-3-partite. Lemma \ref{eq:3.7} implies \( x_{w_1} \geq 1 - \frac{\sqrt{\varepsilon}\,n}{\rho} \) and \( x_{w_2} \geq 1 - \frac{\sqrt{\varepsilon}\,n}{\rho} \).  

Construct \( G' = G - \{ vw : w \in N_{V_2}(u) \} + \{ vw_1, vw_2 \} \). Then \( G' \) is non-3-partite and \( B_{3,k} \)-free. By Lemma \ref{lem2.2}, the spectral radius satisfies:  
$$
\rho(G') - \rho \geq \frac{\boldsymbol{x}^T \big( A(G') - A(G) \big) \boldsymbol{x}}{\boldsymbol{x}^T \boldsymbol{x}} = \frac{2x_v}{\boldsymbol{x}^T \boldsymbol{x}} \left( x_{w_1} + x_{w_2} - \sum_{w \in N_{V_2}(u)} x_w \right).
$$  
Substituting \( x_{w_1}, x_{w_2} \geq 1 - \frac{\sqrt{\varepsilon}\,n}{\rho} \) and \( \sum_{w \in N_{V_2}(u)} x_w \leq |N_{V_2}(u)| \cdot 1 =1 \), we get:  
$$
\rho(G') - \rho \geq \frac{2x_v}{\boldsymbol{x}^T \boldsymbol{x}} \left( 2 - \frac{2\sqrt{\varepsilon}\,n}{\rho} - 1 \right) > 0.
$$  
This contradicts the maximality of \( \rho(G) \). Hence, \( d_{V_i}(u) \geq 2 \) for all \( i \in \{2, 3\} \).  

By analogous reasoning, \( d_{V_i}(v) \geq 2 \) for each \( i \in \{2, 3\} \). Now, for any \( w \in V(G) \setminus \{u, v\} \), \( G - w \) is non-3-partite.  
Without loss of generality, assume \( |V_2| - d_{V_2}(v) \geq |V_3| - d_{V_3}(v) \). Construct  
$$
G' = G \;-\; \left\{ uw : w \in N_G(u) \setminus \{v\} \right\} \;+\; \left\{ uw' : w' \in (V_2 \cup V_3) \setminus N_{V_2}(v) \right\}.
$$  
Then \( G' \) is non-3-partite and \( B_{3,k} \)-free. By Lemma 2.2,  
\[\label{eq:3.11}
\rho(G') - \rho \;\geq\; \frac{\boldsymbol{x}^T \big( A(G') - A(G) \big) \boldsymbol{x}}{\boldsymbol{x}^T \boldsymbol{x}} \;=\; \frac{2x_u}{\boldsymbol{x}^T \boldsymbol{x}} \left( \sum_{\substack{w' \in (V_2 \cup V_3) \setminus N_{V_2}(v)}} x_{w'} \;-\; \sum_{\substack{w \in N_G(u) \setminus \{v\}}} x_w \right). \tag{3.11}
\] 

Since \( |V_2| - d_{V_2}(v) \geq |V_3| - d_{V_3}(v) \), we derive:  
\begin{align}\label{eq:3.12}
|V_2| - d_{V_2}(v) 
&\geq \frac{1}{2} \bigl( |V_2| + |V_3| - d_{V_2}(v) - d_{V_3}(v) \bigr)\notag \\
&= \frac{1}{2} \bigl( n - |V_1| - d(v) + 1 \bigr)\notag \\
&\geq \frac{1}{2} \left( \frac{2n}{3} - 2\sqrt{\varepsilon}\,n - d(v) + 1 \right). \tag{3.12}
\end{align}

For any \( w' \in (V_2 \cup V_3) \setminus N_{V_2}(v) \), \( G - w' \) is non-3-partite, so \( x_{w'} \geq 1 - \frac{\sqrt{\varepsilon}\,n}{\rho} \) (by Lemma 3.7). Combining this with Lemma \ref{lem3.3} and \eqref{eq:3.5}, \eqref{eq:3.7}, \eqref{eq:3.11}, \eqref{eq:3.12}, we obtain  
\begin{align*}
\rho\left( {G}'\right)- \rho  &\geq  \frac{2x_{u}}{\boldsymbol{x}^{T}\boldsymbol{x}}\left[{\left( {\frac{1}{2}\left( {\frac{2n}{3} - 2\sqrt{\varepsilon }n - {d}\left( v\right)  + 1}\right)  + \left( {\frac{1}{3} - 2\sqrt{\varepsilon }}\right) n}\right) \left( {1 - \frac{\sqrt{\varepsilon }n}{\rho}}\right) - {d}\left( u\right)  + 1}\right] \\
&= \frac{2x_{u}}{\boldsymbol{x}^{T}\boldsymbol{x}}\bigg[\left( {\frac{2}{3} - 3\sqrt{\varepsilon}}\right) n-\left( {\frac{2}{3} - 3\sqrt{\varepsilon}}\right)\frac{\sqrt{\varepsilon}n^2}{\rho} - \left( {\frac{1}{2} - \frac{\sqrt{\varepsilon}n}{2\rho}}\right) \left( {{d}\left( u\right)  + {d}\left( v\right) }\right)\\
&\,\,\,\,\,\,\,- \left( {\frac{1}{2} + \frac{\sqrt{\varepsilon}n}{2\rho}}\right) {d}\left( u\right) +O(1)\bigg]\\
&\geq \frac{2x_{u}}{\boldsymbol{x}^{T}\boldsymbol{x}}\bigg[\left( {\frac{2}{3} - 3\sqrt{\varepsilon}}\right) n-\left( {\frac{2}{3} - 3\sqrt{\varepsilon}}\right)\frac{\sqrt{\varepsilon}n^2}{\rho} - \left( {\frac{1}{2} - \frac{\sqrt{\varepsilon}n}{2\rho}}\right)\left(\left( {\frac{2}{3} + 2\sqrt{\varepsilon }}\right) n +2 k\right)\\
&\,\,\,\,\,\,\,-\left( {\frac{1}{2} + \frac{\sqrt{\varepsilon}n}{2\rho}}\right) \left(\left( {\frac{1}{3} + \sqrt{\varepsilon }}\right) n + k\right) +O(1)\bigg]\\
&= \frac{2x_{u}}{\boldsymbol{x}^{T}\boldsymbol{x}}\left[  {\left( {\frac{1}{6} + \frac{\left(7\varepsilon-\sqrt{\varepsilon}\right)n}{2\rho}-\frac{9\sqrt{\varepsilon }}{2}}\right) n + O\left( 1\right) }\right]\\
&> 0,
\end{align*}
a contradiction to the choice of \(G\).

{\bf Case 2.} \(r \geq  4\). If \({d}_{{V}_{i}}\left( u\right)  = 1\) for some \(i \in  \{2,3,\ldots ,r\}\) ({say} \({d}_{{V}_{2}}\left( u\right)= 1\)), then combining with \({N}_{{V}_{2}}\left( u\right)  \cap  {N}_{{V}_{2}}\left( v\right)  \neq \emptyset\) and Claim \ref{c6}, we have \({N}_{{V}_{2}}\left( v\right)  = {V}_{2}\). Now,  by Lemma \ref{lem3.3} and \eqref{eq:3.10},
\begin{align}\notag
\sum_{{i = 3}}^{r}{d}_{{V}_{i}}\left( v\right) = {d}\left( v\right)  - 1 - {d}_{{V}_{2}}\left( v\right)&<\left( {\frac{3}{r} - \frac{3}{{r}^{2}} + 2\sqrt{\varepsilon }}\right) n + {kr} - k + 3 - 1 - \left( {\frac{1}{r} - 2\sqrt{\varepsilon }}\right) n \\ \label{eq:3.13}  
&= \left( {\frac{2}{r} - \frac{3}{{r}^{2}} + 4\sqrt{\varepsilon }}\right) n + {kr} - k + 2.
\end{align}

Without loss of generality, we assume \({d}_{{V}_{3}}\left( v\right) = \min_{3 \leq  i \leq  r}\left\{  {{d}_{{V}_{i}}\left( v\right)}\right\}\). Then by \eqref{eq:3.13},
\begin{align}\notag
{d}_{{V}_{3}}\left( v\right)  \leq  \frac{1}{r - 2}\sum_{{i = 3}}^{r}{d}_{{V}_{i}}\left( v\right)  &\leq  \frac{1}{r - 2}\left[  {\left( {\frac{2}{r} - \frac{3}{{r}^{2}} + 4\sqrt{\varepsilon}}\right) n + {kr} - k + 2}\right] \\ \notag
&\leq  \left( {\frac{1}{r} - \frac{3}{2{r}^{2}} + 2\sqrt{\varepsilon }}\right) n + \frac{{kr} - k + 2}{2},
\end{align}
where the last inequality follows by $r\geq4.$ Therefore, there is at least
\begin{align}\notag
\left| {V}_{3}\right|  - {d}_{{V}_{3}}\left( v\right)  &\geq  \left( {\frac{1}{r} - 2\sqrt{\varepsilon }}\right) n - \left( {\frac{1}{r} - \frac{3}{2{r}^{2}} + 2\sqrt{\varepsilon }}\right) n - \frac{kr - k + 2}{2} \\ \notag
&= \left( {\frac{3}{2{r}^{2}} - 4\sqrt{\varepsilon }}\right) n - \frac{{kr} - k + 2}{2}
\end{align}
vertices in \({V}_{3}\) not adjacent to \(v.\) 

Take \({w}_{1},{w}_{2} \in  {V}_{3} \backslash  {N}_{{V}_{3}}\left( v\right)\). By Claim \ref{c6}, \({w}_{1},{w}_{2}\in N_{V_3}\left( u\right)\), so \(G - {w}_{1}\) and \(G - {w}_{2}\) are non-$r$-partite. Lemma \ref{lem3.7} gives \(x_{{w}_{1}} \geq  1 - \frac{\sqrt{\varepsilon}n}{\rho}\) and \({x}_{{w}_{2}} \geq  1 - \frac{\sqrt{\varepsilon}n}{\rho}\). Construct \({G}' = G - \left\{  {{vw}:w \in  {N}_{{V}_{2}}\left( u\right) }\right\} + \left\{  {{vw_1},v{w}_{2}}\right\}\). Then \({G}'\) is non-$r$-partite and \({B}_{r,k}\)-free. By Lemma \ref{lem2.2},  
$$
\rho(G') - \rho 
\;\geq\; \frac{\boldsymbol{x}^T \big( A(G') - A(G) \big) \boldsymbol{x}}{\boldsymbol{x}^T \boldsymbol{x}} 
\;=\; \frac{2x_v}{\boldsymbol{x}^T \boldsymbol{x}} \left( x_{w_1} + x_{w_2} \;-\; \sum_{\substack{w \in N_{V_2}(u)}} x_w \right).
$$  

Substituting \( x_{w_1}, x_{w_2} \geq 1 - \frac{\sqrt{\varepsilon}\,n}{\rho} \) and \( \sum_{\substack{w \in N_{V_2}(u)}} x_w \leq |N_{V_2}(u)| \cdot 1 =1 \), we get:  
$$
\rho(G') - \rho 
\;\geq\; \frac{2x_v}{\boldsymbol{x}^T \boldsymbol{x}} \left( 2 - \frac{2\sqrt{\varepsilon}\,n}{\rho} - 1 \right) 
\;>\; 0,
$$  
contradicting the maximality of \( \rho(G) \). Hence, \( d_{V_i}(u) \geq 2 \) for all \( i \in \{2, 3, \ldots, r\} \).

By analogous reasoning, \( d_{V_i}(v) \geq 2 \) for each \( i \in \{2, 3, \ldots, r\} \). Now, for any \( w \in V(G) \setminus \{u, v\} \), \( G - w \) is non-\( r \)-partite, so Lemma 3.7 implies \( x_w \geq 1 - \frac{\sqrt{\varepsilon}\,n}{\rho} \).

Since \( r \geq 4 \), by Lemma 3.3 and (3.10):  
$$
d(v) 
\;<\; \left( \frac{3}{r} - \frac{3}{r^2} + 2\sqrt{\varepsilon} \right)n + kr - k + 3 
\;<\; n - \left( \frac{1}{r} + 2\sqrt{\varepsilon} \right)n 
\;\leq\; n - |V_1|.
$$  
This implies \( \left( \bigcup_{i=2}^r V_i \right) \setminus N_G(v) \neq \emptyset \). Without loss of generality, assume \( V_2 \setminus N_{V_2}(v) \neq \emptyset \). Take \( v' \in V_2 \setminus N_{V_2}(v) \), and construct  
$$
G'' = G \;-\; \left\{ uw : w \in N_G(u) \setminus \{v\} \right\} \;+\; \left\{ uw' : w' \in \{v'\} \cup \left( \bigcup_{i=3}^r V_i \right) \right\}.
$$  
Then \( G'' \) is non-\( r \)-partite and \( B_{r,k} \)-free. By Lemmas \ref{lem2.2}, \ref{lem3.3} and \eqref{eq:3.8},
\begin{align*}
\rho \left( {G}'\right)  - \rho &\geq \frac{2x_{u}}{\boldsymbol{x}^{T}\boldsymbol{x}}\left( {\sum_{{{w}' \in  \left\{  {v}'\right\}\cup\left(\bigcup_{{i = 3}}^{r}{V}_{i}\right) }}x_{{w}'} - \sum_{{w \in  {N}_{G}\left( u\right) \backslash  \left\{  v\right\} }}x_{w}}\right) \\
&\geq \frac{2x_{u}}{\boldsymbol{x}^{T}\boldsymbol{x}}\left[{\left( {1 + n - \left| {{V}_{1}}\right|-\left| {V}_{2}\right|}\right) \left( {1 - \frac{\sqrt{\varepsilon}n}{\rho}}\right) - \left( {{d}\left( u\right)  - 1}\right) }\right] \\
&\geq  \frac{2x_{u}}{\boldsymbol{x}^{T}\boldsymbol{x}}\left[ {\left( {\left( {1 - \frac{2}{r} - 4\sqrt{\varepsilon }}\right) n}\right) \left( {1 - \frac{\sqrt{\varepsilon }n}{\rho }}\right) - \left( {\frac{1}{2} - \frac{1}{2r} + \sqrt{\varepsilon }}\right) n +O(1) }\right] \\
&=\frac{2x_{u}}{\boldsymbol{x}^{T}\boldsymbol{x}}\left[ {\left(\frac{1}{2}-\frac{3}{2r}-5\sqrt{\varepsilon}+
\frac{\left(4r\varepsilon+2\sqrt{\varepsilon}-r\sqrt{\varepsilon}\right)n}{r\rho}\right)n+O(1)}\right] \\
&>0,
\end{align*}
a contradiction to the choice of \(G\).

Therefore, there is an index \(i \in  \{ 2,3,\ldots ,r\}\) such that \({N}_{{V}_{i}}\left( u\right)  \cap  {N}_{{V}_{i}}\left( v\right) = \emptyset\).
\end{proof}

Now, we are ready to show Theorem \ref{thm1.5}.
\begin{proof}[\bf Proof of Theorem \ref{thm1.5}.]
By Lemma \ref{lem3.12}, we may assume \( N_{V_2}(u) \cap N_{V_2}(v) = \emptyset \). Without loss of generality, let \( x_u \leq x_v \). Denote \( d_{V_2}(u) = s \); we claim \( s = 1 \).  

\noindent{\bf Justification for \( s = 1 \):}   
Since \( V(G) = V_1 \cup V_2 \cup \cdots \cup V_r \) maximizes \( \sum_{1 \leq i < j \leq r} e(V_i, V_j) \), we have \( s \geq d_{V_1}(u) \geq 1 \). If \( s \geq 2 \), take \( w' \in N_{V_2}(u) \) and construct  
$
G' = G - uw' + vw'.
$ 
Then \( G' \) is non-\( r \)-partite and \( B_{r,k} \)-free. Since \( x_u \leq x_v \), Lemma \ref{lem2.6} implies \( \rho(G') > \rho(G) \), contradicting the maximality of \( \rho(G) \). Thus, \( s = 1 \).

\noindent{\bf Structure of \( G \):} \ Now, \( G \) is a subgraph of \( K_{|V_1|, |V_2|, \ldots, |V_r|}^1 \), where \( K_{|V_1|, |V_2|, \ldots, |V_r|}^1 \) is constructed as follows:  
\begin{enumerate}
  \item Start with a complete \( r \)-partite graph with partitions \( V_1, V_2, \ldots, V_r \).
  \item Add an edge between \( u, v \in V_1 \). 
  \item Delete \( |V_2| - 1 \) edges between \( u \) and \( V_2 \), and 1 edge between \( v \) and \( V_2 \), such that \( u \) and \( v \) have no common neighbors in \( V_2 \). 
\end{enumerate}

\noindent{\bf Spectral Radius Analysis:}\ Clearly, \( K_{|V_1|, |V_2|, \ldots, |V_r|}^1 \) is non-\( r \)-partite and \( B_{r,k} \)-free. By the Perron-Frobenius theorem:  
$$
\rho(G) \leq \rho\left( K_{|V_1|, |V_2|, \ldots, |V_r|}^1 \right),
$$ 
with equality if and only if \( G = K_{|V_1|, |V_2|, \ldots, |V_r|}^1 \).  Further, by Lemma \ref{lem2.7}:  
$$
\rho(G) \leq \rho\left( Y_r(n) \right),
$$ 
with equality if and only if \( G = Y_r(n) \). 
\end{proof}
\section{\normalsize Further discussions}

For any fixed graph \( H \), exploring the relationship between \( \text{SPEX}(n, H) \) and \( \text{EX}(n, H) \) constitutes an intriguing research problem within both Tur\'an-type problems and spectral Tur\'an-type problems. Wang, Kang, and Xue \cite{WKX2023} demonstrated that for \( r \geq 2 \) and sufficiently large \( n \), if \( H \) is a graph satisfying \( \text{ex}(n, H) = e(T_r(n)) + O(1) \), then \( \text{SPEX}(n, H) \subseteq \text{EX}(n, H) \). This result resolved a conjecture proposed by Cioab\u{a}, Desai, and Tait \cite{CDT2022}. In 2023, Liu and Ning formulated the following problem:  

\begin{pb}[\cite{LN2023}]\label{pb4.1} 
Let \( H \) be an arbitrary graph. Characterize all graphs \( H \) for which \( \text{SPEX}(n, H) \subseteq \text{EX}(n, H) \).  
\end{pb}

Recently, Fang, Tait, and Zhai \cite{FTZ2025} solved Problem \ref{pb4.1} for graphs \( H \) with \( \chi(H) = r + 1 \geq 3 \) and \( \text{ex}(n, H) < e(T_r(n)) + \left\lfloor \frac{2r}{n} \right\rfloor \), thereby extending the work of Wang, Kang, and Xue \cite{WKX2023}.
Amin et al. \cite{AFGS} established the value of $\text{EX}_{r + 1}(n, K_{r + 1})$ for $r \geq 2$. Subsequently, Lin, Ning, and Wu \cite{LNW}, as well as Li and Peng \cite{LP2023}, determined $\text{SPEX}_{r + 1}(n, K_{r + 1})$ for $r \geq 2$. Through a comparison of $\text{EX}_{r + 1}(n, K_{r + 1})$ and $\text{SPEX}_{r + 1}(n, K_{r + 1})$, we observe that $\text{SPEX}_{r + 1}(n, K_{r + 1}) \subseteq \text{EX}_{r + 1}(n, K_{r + 1})$ when $r \geq 2$. Ren et al. \cite{RWWY} identified $\text{EX}_3(n, C_{2k + 1})$ for $k \geq 2$, while Zhang and Zhao \cite{ZZ2023} specified $\text{SPEX}_3(n, C_{2k + 1})$ for $k \geq 2$. Analogously, by comparing $\text{EX}_3(n, C_{2k + 1})$ and $\text{SPEX}_3(n, C_{2k + 1})$, we find that $\text{SPEX}_3(n, C_{2k + 1}) \subseteq \text{EX}_3(n, C_{2k + 1})$ for $k \geq 2$. For additional findings in this research vein, one may consult our recent papers \cite{CLL2025,HL2024,LL2024}.   

In our previous work \cite{YL}, we determined \(\text{EX}_{r + 1}(n, B_{r,k})\) for \(r \geq 3\), \(k \geq 1\), and sufficiently large \(n\). In this paper, we consider the spectral version of the result from \cite{YL} and determine \(\text{SPEX}_{r + 1}(n, B_{r,k})\) for \(r \geq 3\), \(k \geq 1\), and sufficiently large \(n\). By comparing \(\text{EX}_{r + 1}(n, B_{r,k})\) and \(\text{SPEX}_{r + 1}(n, B_{r,k})\), we find that the phenomenon \(\text{SPEX}_{r + 1}(n, B_{r,k}) \subseteq \text{EX}_{r + 1}(n, B_{r,k})\) persists for all \(r \geq 3\), \(k \geq 1\), and sufficiently large \(n\). Motivated by Problem 3 and the aforementioned facts, we propose the following problem. 

\begin{pb}\label{pb4.2}
Is it true that for all integers \(r \geq 2\) and all color-critical graphs \(H\) with \(\chi(H)=r + 1\), we always have \(\text{SPEX}_{r + 1}(n, H) \subseteq \text{EX}_{r + 1}(n, H)\) for sufficiently large \(n?\)
\end{pb}

Another issue drawing attention is a classical result by Nosal \cite{EN1970} (see, e.g., \cite{VN2002}). It states that for a triangle-free graph \( G \) with \( m \) edges, the spectral radius satisfies \( \rho(G) \leq \sqrt{m} \).  In the literature, Nosal’s theorem has seen various generalizations to scenarios where other subgraphs \( F \) are forbidden. For those interested in exploring this line of research further, a recent paper \cite{LZZ2025} serves as a valuable reference. 

In 2001, Nikiforov \cite{VN2021} obtained the result below:
\begin{thm}[\cite{VN2021}]
If $G$ is an $m$-edge and $B_{2,k}$-free graph, then $\rho(G)\leqslant \sqrt{m},$ with equality if and only if $G$ is a complete bipartite graph.
\end{thm}

Nikiforov’s theorem above inspires the investigation of the maximum spectral radius of non-bipartite \( B_{2,k} \)-free graphs with \( m \) edges (i.e., of size \( m \)), along with characterizing the corresponding spectral extremal graphs. Building on this insight, Liu and Miao \cite{LM2025} proposed the following conjecture. 
\begin{conj}
Let \( G \) be a non-bipartite \( B_{2,k} \)-free graph \( (k\geqslant 2) \) with \( m \) edges. Then, the spectral radius satisfies \( \rho(G)\leq\rho(U) \). Equality holds if and only if \( G \cong U \), where \( U \) is constructed by attaching \( m - 3 \) pendant edges to precisely one vertex of the $3$-cycle \( C_3 \). 
\end{conj}

We also notice that Li, Liu and Feng \cite{LLF2022} put forward the following conjecture for the generalized book graph $B_{r,k}$:
\begin{conj}[\cite{LLF2022}]\label{cj2}
Let \( r \geq 2 \), \( k \geq 1 \) be fixed constants, and let \( m \) be sufficiently large. If \( G \) is a \( B_{r,k} \)-free graph with \( m \) edges, then \( \rho(G) \leq \sqrt{\left(1 - \frac{1}{r}\right) \cdot 2m}\). 
\end{conj}
Conjecture \ref{cj2} is interesting in its own right, and there are additional aspects of it worthy exploring. In fact, we do not know whether the upper bound in Conjecture \ref{cj2} is tight, nor is the structure of the corresponding extremal graph clear. 

Note that both the book graph \( B_{2,k} \) and the generalized book graph \( B_{r,k} \) are color-critical. In view of the above mathematical phenomenon, we can consider the following general problem.
\begin{pb}
Given a color-critical graph \( H \) with \( \chi(H)=r + 1 \), determine the graphs that attain the maximum spectral radius among all \( H \)-free graphs with \( m \) edges.
\end{pb}
We intend to explore these problems in the near future. 

\section*{\normalsize Declaration of competing interest}
The authors declare that they have no known competing financial interests or personal relationships that could have influenced the work reported in this paper.

\section*{\normalsize Acknowledgement}
S.L. is financially supported by the National Natural Science Foundation of China (Grant Nos. 12171190, 11671164), the Special Fund for Basic Scientific Research of Central Colleges (Grant No. CCNU25JC006) and the Open Research Fund of Key Laboratory of Nonlinear Analysis \& Applications (CCNU), Ministry of Education of China (Grant No. NAA2025ORG010). 

We also thank Dr. Bing Wang \cite{wang}, who noticed our preprint \cite{YL} and shared with us an unpublished manuscript from his team. Upon reviewing their work, we confirmed that they independently obtained the main result (Theorem~\ref{thm1.5}) of this paper as well.
\section*{\normalsize Data availability}
No data were utilized for the research that is described within this article.

\end{document}